\newif\ifjournal
\journalfalse

\newif\ifray
\rayfalse

\newif\ifdraft
\draftfalse

\ifjournal
\ifdraft
  \documentclass[times,doublespace]{nlaauth}
\else
  \documentclass[times]{nlaauth}
\fi

\else
\documentclass[times]{article}

\usepackage{hyperref}

\newcommand{\tabsize}{\fontsize{9}{9.5pt}\selectfont}

\fi % \ifjournal

\usepackage{algpseudocode}
\usepackage{algorithm}
\usepackage{amsmath,amssymb}
\usepackage{bm}
\usepackage{booktabs}
\usepackage{multicol}
\usepackage{multirow}

\ifdraft
\usepackage[pagewise,modulo]{lineno}
\linenumbers
\fi

\DeclareMathOperator*{\argmin}{arg\,min}
\DeclareMathOperator*{\diag}{diag}
\DeclareMathOperator{\lump}{lump}

\usepackage[pdftex,dvipsnames]{xcolor}
\ifray
\else
\usepackage{pgfplots}
\pgfdeclarelayer{bg}    % declare background layer
\pgfsetlayers{bg,main}  % set the order of the layers (main is the standard layer)
\fi

\usepackage[colorinlistoftodos,prependcaption,textsize=small]{todonotes}

\usepackage{xargs}
\usepackage{soul}
\newcommandx{\fix}     [3][1=]{\todo[linecolor=RedOrange,backgroundcolor=RedOrange!25,bordercolor=RedOrange,#1]{\textbf{#2: }#3}}
\newcommandx{\unsure}  [3][1=]{\todo[linecolor=Cerulean,backgroundcolor=Cerulean!25,bordercolor=Cerulean,#1]{\textbf{#2: }#3}}
\newcommandx{\improve} [3][1=]{\todo[linecolor=OliveGreen,backgroundcolor=OliveGreen!25,bordercolor=OliveGreen,#1]{\textbf{#2: }#3}}
\newcommandx{\info}    [3][1=]{\todo[linecolor=Plum,backgroundcolor=Plum!25,bordercolor=Plum,#1]{\textbf{#2: }#3}}
\newcommandx{\RT}[1]{#1}
\newcommandx{\changehl}[2]    {\texthl{#1}\change{#2}}
\ifjournal
\newcommand{\qq}{\ensuremath{\bm{Q}_2\,}--\ensuremath{\,\bm{Q}_1\,}}
\else
\newcommand{\qq}{\ensuremath{Q_2\,}--\ensuremath{\,Q_1\,}}
\fi

\newcommand{\func}[1]{\texttt{#1}}

\begin{document}

\ifjournal
\runningheads{~A.~Prokopenko,~R.~S.~Tuminaro}{AMG for \qq
discretization}
\fi

\title{An algebraic multigrid method for \qq mixed discretizations of
the Navier-Stokes equations}

\ifjournal
\author{A.~Prokopenko\affil{1},
R.~Tuminaro\affil{2}\comma\corrauth}

\address{
\affilnum{1}\ Center for Computing Research, Sandia National Laboratories,
Albuquerque, NM 87185 \break
\affilnum{2}\ Center for Computing Research, Sandia National Laboratories,
Livermore, CA 94551
}

\corraddr{Sandia National Laboratories,
P.O. Box 969, MS 9159, Livermore, CA 94551, USA ({\tt rstumin@sandia.gov}).}

\keywords{preconditioning; algebraic multigrid; Navier-Stokes equations; mixed
finite element discretizations}

\else % \ifjournal

\author{Andrey Prokopenko\thanks{Center for Computing Research, Sandia National Laboratories, Albuquerque, NM 87185
  (\texttt{aprokop@sandia.gov})}\; and
Raymond S. Tuminaro\thanks{Center for Computing Research, Sandia National Laboratories, Livermore, CA 94551
(\texttt{rstumin@sandia.gov})}}

\date{}

\fi % \ifjournal

\ifjournal
\else % \ifjournal
% For preprint, put abstract after \maketitle

\maketitle

\begin{abstract}
Algebraic multigrid (AMG) preconditioners are considered for
discretized systems of partial differential equations (PDEs) where unknowns
associated with different physical quantities are not necessarily co-located
at mesh points.  Specifically, we investigate a \qq mixed finite element
discretization of the incompressible Navier-Stokes equations where the
number of velocity nodes is much greater than the number of pressure
nodes. Consequently, some velocity degrees-of-freedom (dofs) are
defined at spatial locations where there are no corresponding
pressure dofs.  Thus, AMG approaches leveraging this co-located
structure are not applicable.
This paper instead proposes an automatic AMG coarsening that mimics certain
pressure/velocity dof relationships of the \qq discretization.
The main idea is to first automatically define coarse pressures in a
somewhat standard AMG fashion and then to carefully (but automatically) choose coarse
velocity unknowns so that the spatial location relationship between pressure
and velocity dofs resembles that on the finest grid.  To
define coefficients
within the inter-grid transfers, an energy minimization AMG (EMIN-AMG)
is utilized. EMIN-AMG is not tied to specific coarsening schemes and
grid transfer sparsity patterns, and so it is
applicable to the proposed coarsening.
Numerical results highlighting solver performance are given on Stokes and
incompressible Navier-Stokes problems.

% as well as on a potential
% formulation of a resistive magnetohydrodynamics (MHD) system.

% In the MHD
% example, we use a family of recently proposed multigrid smoothers that
% correspond to
% extensions of commonly used smoothers for incompressible flow.
 % (Vanka smoothing and Braess-Sarazin smoothing).
% On structured grids
% the AMG scheme more closely resembles geometric multigrid

\end{abstract}

\fi % \ifjournal

\footnotetext[1]{Sandia is a multiprogram laboratory operated by Sandia Corporation, a Lockheed Martin Company, for the United States Department of Energy under contract DE-AC04-94-AL85000.  Part of this material is based upon work supported by the U.S. Department of Energy, Office of Science, Office of Advanced Scientific Computing Research, Applied Mathematics program.}

\section{Introduction}\label{s:intro}
Multigrid methods are among the most efficient algorithms for solving sparse
linear systems arising from discretized elliptic partial differential
equations (PDE)s~\cite{BrHeMc2000,TrOoSc2001}.
Rapid convergence requires that the algorithm's relaxation phase
complement the coarse grid correction phase.
Relaxation focuses on reducing oscillatory error components, often
via a simple iteration.
%The coarse problem complements relaxation by reducing smooth error components.
The coarse phase projects %problem can be formulated by projecting
a residual equation to a coarse space and interpolates an associated
coarse {\it solution} to correct the
approximation.  A hierarchy  of resolutions result when the
coarse solution is approximated by a recursive multigrid invocation.
This paper is concerned with mixed finite element
discretizations of PDE systems and an algebraic multigrid (AMG) method, which
automatically constructs a mesh hierarchy and grid transfers.
%Algebraic multigrid is desirable for unstructured mesh applications as
%it does not require detailed mesh connectivity data.  However,
For PDE systems, constructing AMG components with the desired
%relaxation/coarse correction
complementary properties is
challenging, especially when there
is strong coupling between different types of unknowns (e.g., pressures
and velocities). %Unlike many discrete representations of scalar
%elliptic PDEs,
Matrices that result from PDE systems are frequently far
from the M-matrices that are generally more amenable to standard AMG methods.
\RT{Among applications involving PDE systems, the saddle point nature of incompressible
flow problems introduces additional complications. Standard multigrid relaxation algorithms do not typically smooth errors appropriately and algebraic multigrid procedures for generating grid transfers that partially rely on positive-definite matrix properties (such as smoothed aggregation) have difficulties as well.}
Further, traditional algebraic multigrid methods
do not distinguish between unknown types.  Thus, they might
produce odd interpolation operators with stencils that mix unknown types,
e.g., a fine pressure might interpolate from a set of coarse velocities.

Geometric multigrid methods have been developed for a range of PDE
systems and discretizations~(e.g., \cite{TrOoSc2001}). Much of this work has centered on
effective relaxation techniques for specific PDE systems and so ideally
one would like to leverage these techniques within an AMG
approach.  In geometric multigrid, inter-grid transfers are
based on a geometrical relationship between meshes, such as using linear
interpolation to transfer solutions across meshes.
AMG methods are an attractive alternative for applications with complex meshes
or features as formulating a mesh hierarchy and geometric inter-grid transfers
is sometimes challenging.

One AMG approach (referred to as the point-based approach in \cite{FuSt02})
for
co-located PDE systems centers on coarsening spatial
locations.
%In most PDE cases where multigrid has been applied to an entire PDE
%system, the unknowns associated with the matrix equation have either
%been co-located or for more general discretizations geometric multigrid has
%been more prevalent.
% been employed or when algebraic multigrid has been used
%Most earlier developments
%centered on geometric multigrid.  Generally, multigrid hierarchies and
%grid transfers can be defined
%in fairly natural ways when the fine mesh is obtained by uniformly
%refining a coarse mesh. For example, when finite elements are used,
%it is common to define inter-grid transfers by effectively evaluating
%coarse grid basis functions at fine grid points. That is, finite
%element spaces naturally define the grid transfers.
%Finally, one may still need to design an appropriate smoother.
%It should, however,
%be noted that such grid transfers might not be appropriate
%(e.g.  when coarse and fine meshes that are not nested,
%non-finite element applications, anisotropic problems,
%applications with interfaces or large material variations, and non-elliptic
%situations such as Helmholz systems) and this has been the subject
%of a significant body of geometric multigrid research.
%As an alternative to geometric multigrid, there is a significant body
%of algebraic multigrid research associated with PDE systems having
%co-located unknowns. % centering on block matrices.
%The algebraic multigrid situation is generally more complex.  When unknowns
%are co-located, a natural AMG approach centers on block matrices.
By co-located, we refer to discretizations where degrees-of-freedom (dofs)
associated with different physical quantities (e.g., pressures and
 electric fields) are situated at the same spatial location.
For a co-located representation of the incompressible Navier-Stokes equations,
each velocity dof is
associated with a spatial location where there is a corresponding pressure dof.
The discretization matrix
can be viewed as a block operator with constant sized blocks, each
corresponding to a unique spatial location. Standard AMG coarsening
constructs a graph where each vertex is associated with one block matrix
row and weighted edges between vertex $i$ and $j$ are
defined based on some norm (or quantity) corresponding to the sub-matrix
in the $i$th block row and the $j$th block column. Standard
coarsening algorithms can then be employed accompanied by schemes
(typically modest modifications of AMG methods for scalar PDEs)
to define coefficients for the grid transfer operators such that
coarse level dofs retain the co-located structure.
Unfortunately, this PDE system AMG approach is not
possible for mixed finite element discretizations in which fine
level dofs are not co-located. %\textit{e.g.} \qq discretization.

A second AMG approach
(referred to as the unknown-based approach in \cite{FuSt02})
also considers a block matrix
representation, but now blocks are defined by different physical
quantities (e.g., pressures)
and different equation types (e.g., Navier-Stokes momentum equations,
incompressibility condition, Maxwell-Faraday equations).
Specifically, the matrix can be written as
\begin{equation}
\begin{bmatrix}
A_{11} & A_{12} & \cdots & A_{1k} \\
A_{21} & A_{22} & \cdots & A_{2k} \\
\vdots & \vdots & \ddots & \vdots \\
A_{k1} & A_{k2} & \cdots & A_{kk} \\
\end{bmatrix}
\begin{bmatrix}
u_1 \\
u_2 \\
\vdots \\
u_k \\
\end{bmatrix}
=
\begin{bmatrix}
f_1 \\
f_2 \\
\vdots \\
f_k \\
\end{bmatrix},
\label{eq:ns_systxm}
\end{equation}
where each $u_i$ corresponds to a field
(e.g.,
%In the case of the incompressible Navier-Stokes equations,
%$k$, might be four and the $f_i$'s would correspond
velocities in the $x, y,$ and $z$ directions).
% and pressure for the incompressible Navier-Stokes equations).

An AMG algorithm can be developed by requiring
that grid transfers have a block diagonal form:
\begin{equation}
  P =
  \begin{pmatrix}
P_1 &   &   &   &    \\
    & . &   &   &    \\
    &   & . &   &    \\
    &   &   & . &    \\
    &   &   &   & P_k
\end{pmatrix}.
\label{eq:general_block_p}
\end{equation}
Thus, coarse unknowns of one type do not directly
influence interpolated unknowns of another type. % (e.g., velocities).
%It is important to note, however, that the entire PDE system
%is projected (using \eqref{eq:general_block_p}),
%and so coarse problems still retain off-diagonal coupling.
The $P_i$ can be produced by $k$ separate invocations of an AMG
method applied to matrices $\widetilde{A}_{ii}$. When $\widetilde{A}_{ii} =
A_{ii}$, the AMG grid transfers ignore the PDE cross-coupling
and so this might be problematic when the coupling
is relatively strong.  For mixed finite element representations
that include an incompressibility condition,
there is often a diagonal block of the PDE system that is identically
zero, and so $\widetilde{A}_{ii}$ must be defined in another
way (e.g., using Schur complement ideas).

This article follows a block interpolation approach.
%where AMG is applied to different $\widetilde{A}_{ii}$ operators.
Our main innovation is to couple AMG invocations in a limited way.
%Instead, the AMG coarsening algorithms are linked.
This is done to {\it mimic} certain features of a \qq discretization
on coarse grids as %.  It is well known that
\qq discretizations of the
incompressible Navier-Stokes equations satisfy inf-sup conditions
needed to produce stable
discretizations~\cite{brezzi91,Gunzburger89,Elman2005}.
It is generally desirable that all discretization operators within
a multigrid method be stable
as unstable coarse operators typically pose multigrid convergence problems
(e.g., relaxation schemes may diverge). Obviously, geometric multigrid
algorithms  employing \qq discretizations on all levels satisfy
inf-sup conditions throughout the hierarchy.
We seek to emulate this within an AMG method by selecting coarse points
in a fashion that loosely resembles coarse points within
a geometric multigrid method.
In particular, the selection of coarse velocities
(used to define velocity interpolation) is obtained by first including
velocities co-located with coarse pressures determined during a
prior pressure AMG invocation. % (to define pressure interpolation).
The set of coarse
velocity unknowns is then augmented by velocity unknowns
located at approximate {\it mid-points} between the coarse pressure
unknowns.
%While the overall approach is not completely algebraic in that it
%employs coordinate locations, it requires
%users to only provide minimal (coordinate location) geometric
%information.
Numerical results will be given to demonstrate the effectiveness of this
strategy when it is coupled to an energy minimization AMG
algorithm~\cite{Olson2011}, though the coarsening does not guarantee
that resulting coarse discretizations satisfy inf-sup conditions.
One key idea is that it is relatively
easy to employ fairly general coarsening and grid transfer patterns
within an energy minimization AMG (EMIN-AMG) framework. This is because this
AMG variant (unlike most others) is fairly flexible with respect to different
coarsenings and sparsity patterns. Other methods such as bootstrap AMG also
have this property~\cite{bootstrap}.
In addition to addressing stability concerns, it is expected that a
careful choice of coarse variables will enhance the use of
physics-based relaxation methods that rely heavily on sub-matrices
accurately mirroring properties of associated PDE
%Additionally, we note that many specialized methods that can be used
%as relaxation techniques heavily rely on sub-matrix
%properties mirroring those of the corresponding PDE
operators~\cite{Elman2005,Wesseling2001}.

The correlation of coarse unknowns between separate prolongator operators
%(i.e., $P_i$ in \eqref{eq:general_block_p})
has been considered in different settings.
In \cite{wabro2004}, multigrid transfers are motivated by recognizing
that a standard geometric coarsening of the velocity mesh associated
with a $\bm{P}_1$iso$\bm{P}_2\,$--$\,\bm{P}_1$ discretization corresponds
to the pressure mesh. This leads to a {\it shift} strategy whereby
the first coarsening of the velocity mesh is given by the
pressure mesh (and the first velocity grid transfer could simply inject
velocities located at nodes of the pressure mesh). Subsequent velocity
grid transfers to the $\ell$th multigrid level use the AMG generated
pressure grid transfer operator to the $\ell\hskip -.04in -\hskip -.04in1$th
level. Coarse level stability is considered and alternative
strategies retaining more fine level velocities are also proposed
for \qq discretizations. % with different cost/convergence tradeoffs.
While good convergence rates are obtained on several non-trivial problems,
the approach uses a somewhat restrictive coarsening procedure that
can result in multigrid cycles with an expensive cost. Our proposed
algorithm seeks increased flexibility where one can use more aggressive coarsening
rates to reduce the number of levels and the number of nonzeros per row
(i.e., fill in) of the coarse operators, particularly important in a parallel
setting.

Multigrid methods based on element agglomeration (AMGe) have also been
considered for mixed finite element problems and mimetic discretizations~\cite{amgeone,amgetwo,wabroamge}.
In AMGe methods, additional topological information is maintained
throughout all multigrid levels. While initially providing and
maintaining this topological information has some challenges, its presence
on coarse hierarchy levels can be leveraged in developing mixed fixed element
AMG methods that satisfy key topological properties. Additionally,  special multigrid methods for solid mechanics
problems with contact and slide surfaces bear some resemblance to our
proposed %\qq finite element
solver in that the coarsening of constraints
associated with contact are correlated to the coarsening of
displacement dofs~\cite{adams_uzawa}.
As another example, AMG and fluid-structure interactions are described
in \cite{gee,Yang20115367}.

%In the \qq case, each pressure unknown
%is co-located with velocity unknowns. Additionally, velocity unknowns
%are defined at spatial locations associated with mid-points between
%adjacent pressure unknowns. This particular arrangement of unknowns
%results from the \qq finite elements, and these elements are generally chosen
%because of their discretization stability properties, \textit{i.e.} they satisfy inf-sup
%conditions. In mirroring this structure, there is an expectation, though
%definitely not a guarantee,  that coarse level discretizations will also be
%stable.  We accomplish this goal by first coarsening pressures and then modifying
%the coarsening of velocities so that coarse pressures have co-located
%coarse velocities and that additional coarse velocities reside at {\it near}
%mid-point locations between coarse pressures.

Before describing our AMG approach, we note that physics-based
preconditioners are an alternative for PDE systems.
These techniques also consider a block system
such as \eqref{eq:ns_systxm}
%, but avoid some challenges by instead
%only applying multigrid to sub-matrices associated with just velocities
%and/or with just pressures.
and can be viewed as
approximate block factorizations %of \eqref{eq:ns_systxm}
involving
Schur complement approximations.
%based on PDE properties of individual operators.
Here, AMG is used to approximate sub-matrix inverses needed
within the approximate block factors.
%There are many connections between these methods and some specialized PDE
%smoothing techniques developed for geometric multigrid.
%As AMG is not applied to the entire
%system, challenges associated with mixed finite elements are avoided.
While physics-based preconditioners are effective
and scalable, their convergence behavior is tied
to the Schur complement approximations.
The current paper is motivated by %our observation of %that there are
situations where monolithic %AMG %(applied to discretizations resulting
multigrid
%in co-located unknowns)
can outperform %require noticeably less run time than
several approximate block factorization preconditioners, though
cases also exists where approximate block factorization preconditioners
are better~\RT{\cite{CyShTuPaCh2012}}.

%%As with geometric multigrid methods, appropriate relaxation techniques must then be utilized.
%%For example, the Gauss-Seidel method may no longer be
%%an effective smoother of errors and may not be even defined when the
%%matrix diagonal contains zero entries (as would be the case for a
%%\qq discretization of the incompressible Navier-Stokes equations).
%%%
%%Smoothers specific to particular PDE systems have been developed, often
%%relying on PDE properties of sub-matrices within the overall system.
%%For example, Braess-Sarazin relaxation views the discrete incompressible
%%Navier-Stokes operator as a $2 \times 2$ block matrix where the four
%%sub-matrices represent different degree-of-freedom(dof)/equation couplings:
%%velocity dofs/momentum equations, pressure dofs/momentum equations,
%%velocity dofs/incompressibility equation, and pressure
%%dofs/incompressibility equation (when explicit pressure stabilization
%%is employed). A key aspect of the smoother is a Schur complement approximation
%%that is based on the properties of the individual operators.
%%

%In Section~\ref{s:discretization}, we describe the \qq discretization for
%several problem formulations. In Section~\ref{s:algorithm}, we introduce the
%energy-minimizing multigrid framework, and describe the algorithms that are used
%to construct the sparsity patterns.  In Section~\ref{s:results}, we present
%numerical evidence for the algorithms. Finally, in Section~\ref{s:conclusion},
%concluding remarks are made.

The paper is organized as follows. Section~\ref{s:discretization} briefly
summarizes the considered PDEs and the \qq discretization.
Section~\ref{s:example} illustrates a potential multigrid stability pitfall when
velocity and pressure coarsenings are not correlated. Section~\ref{s:algorithm}
describes a new algorithm for constructing grid transfers while
Section~\ref{s:smoothers} completes the AMG description by detailing the
multigrid relaxation. Numerical results are given in Section~\ref{s:results}
followed by the conclusions in Section~\ref{s:conclusion}.
% Finally, appendices are given that include details of the graph heuristics.

While our focus is on
\qq discretizations, it is hoped that the techniques in this paper could be
extended to other situations where one wishes to preserve specific properties of
the fine level PDE discretization (e.g., mimetic representations) throughout a
multigrid hierarchy.
%While this particular paper is focused on \qq discretizations,
%there are many situations where PDE discretizations are desirable
%even though they do not result in co-located unknowns. This is often
%motivated by discretization stability, accuracy considerations, or
%by a requirement that the discrete system preserve key properties,
%e.g. conservation of certain quantities or mimetic discretizations
%having discrete structures that mimic vector calculus identities.

%In this specific case, the number of velocity nodes is much
%greater than the number of pressure nodes. Consequently, some velocity
%dofs are associated with spatial locations where there
%are no corresponding pressure dofs.  Thus, standard coarsening is no longer
%natural as the matrix cannot be viewed as a constant sized block operator
%with blocks defined by spatial location.

\section{Problem formulation}\label{s:discretization}
\subsection{Discretization of the Stokes equations}
The Stokes equations describe a model of viscous flow. Given a two- or
three-dimensional domain $\Omega$, let $\bm{u}$ be a vector-valued function
representing velocity, $p$ be a scalar function representing pressure, \RT{and
$\bm{f}$ be a forcing term}. The Stokes system is written as
\begin{equation}
  \begin{array}{rcl}
    -\nabla^2 \bm{u} + \nabla p &=& \bm{f} \quad\mbox{in}\;\Omega, \\
    \nabla \cdot \bm{u} &=& 0 \quad\mbox{in}\;\Omega.
  \end{array}
\label{eq:stokes_system}
\end{equation}
It is complemented by a set of boundary conditions on $\partial\Omega =
\partial\Omega_D \cup \partial\Omega_N$ given by
\begin{equation}
  \bm{u} = \bm{w} \;\mbox{on}\; \partial\Omega_D, \qquad
  \frac{\partial\bm{u}}{\partial\bm{n}} - \bm{n}p = \bm{s}
  \;\mbox{on}\; \partial\Omega_N,
\label{eq:stokes_bc}
\end{equation}
where $\bm{n}$ is the outward-pointing normal to the boundary and both $\bm{w}$
and $\bm{s}$ are given. If the
velocity is specified on the whole boundary, i.e. $\partial\Omega_D =
\partial\Omega$, then the compatibility condition
\begin{equation}
  \int\limits_{\partial\Omega} \bm{w} \cdot \bm{n} = 0
\end{equation}
must be satisfied. In this case, pressure is unique only up to a constant.
%sometimes referred as the \textit{hydrostatic} pressure.

A weak formulation of the Stokes equations is written as:

Find $\bm{u} \in \mathbf{H}^1(\Omega)$ and $p \in L_2(\Omega)$ such that
\begin{equation}
  \begin{array}{rcll}
  \displaystyle\int\limits_{\Omega} \nabla \bm{u} : \nabla \bm{v} - \int\limits_{\Omega} p \nabla
    \cdot \bm{v} &=& \displaystyle\int\limits_{\partial\Omega_N} \bm{s} \cdot
    \bm{v}  + \displaystyle\int\limits_{\Omega} \bm{f} \cdot \bm{v} &
  \mbox{for all}\; \bm{v} \in \mathbf{H}_0^1(\Omega), \\
  \displaystyle\int\limits_{\Omega} q\nabla\cdot\bm{u} &=& 0 & \mbox{for all}\; q \in
    L_2(\Omega),
  \end{array}
\label{eq:stokes_weak_formulation}
\end{equation}
where
\begin{equation}
  \begin{array}{lcl}
    \mathbf{H}^1(\Omega) &=& \left\{\bm{u} \in \mathcal{H}^1(\Omega)^d \;|\; \bm{u} =
      \bm{w} \;\mbox{on}\; \partial\Omega_D\right\}, \\
    \mathbf{H}_0^1(\Omega) &=& \left\{\bm{v} \in \mathcal{H}^1(\Omega)^d \;|\;
      \bm{v} = \bm{0} \;\mbox{on}\; \partial\Omega_D\right\},
  \end{array}
\end{equation}
with $d = 2$ or $d = 3$ being the spatial dimension.

The weak formulation~\eqref{eq:stokes_weak_formulation} is discretized using
finite-dimensional subspaces $\mathbf{X}^h \subset \mathbf{H}^1(\Omega)$,
$\mathbf{X}_0^h \subset \mathbf{H}_0^1 (\Omega)$ and $M^h \subset L_2(\Omega)$,
and is written as:

Find $\mathbf{u}_h \in \mathbf{X}^h$ and $p_h \in M^h$ such that
\begin{equation}
  \begin{array}{rcll}
  \displaystyle\int\limits_{\Omega} \nabla \mathbf{u}_h : \nabla \mathbf{v}_h -
  \displaystyle\int\limits_{\Omega} p_h \nabla
    \cdot \mathbf{v}_h &=& \displaystyle\int\limits_{\partial\Omega_N} \mathbf{s} \cdot
  \mathbf{v}_h + \displaystyle\int\limits_{\Omega} \bm{f} \cdot \mathbf{v}_h &
    \mbox{for all}\; \mathbf{v}_h \in \mathbf{X}_0^h, \\
  \displaystyle\int\limits_{\Omega} q_h\nabla\cdot\mathbf{u}_h &=& 0 & \mbox{for all}\; q_h \in M^h.
  \end{array}
\label{eq:stokes_discrete_formulation}
\end{equation}
Assuming that all $d$ components of velocity are approximated using \RT{the same
scalar finite element space}, this leads to a discrete Stokes problem
\begin{equation}
\begin{pmatrix} A_S & B^T \\ B & O \end{pmatrix}
  \begin{pmatrix} \mathbf{u} \\ \mathbf{p} \end{pmatrix} =
  \begin{pmatrix} \mathbf{f}_{\mathbf{u}} \\ \mathbf{f}_{\mathbf{p}} \end{pmatrix}.
\label{eq:stokes_block_system}
\end{equation}
Matrices $A_S$ and $B$ are block matrices, and have the
following form for $d = 2$:

\RT{$A_S = \begin{pmatrix} \widehat{A}_S & O \\ O & \widehat{A}_S \end{pmatrix}$},
$B   = \begin{pmatrix} B_x & B_y \end{pmatrix}$ and
\begin{equation}
  \begin{array}{cclccl}
    \widehat{A}_S &=& [a_{ij}],   & a_{ij}   &=& \displaystyle\int\limits_{\Omega} \nabla
    \phi_i \cdot \nabla \phi_j, \\
    B_x &=& [b_{x,ki}], & b_{x,ki} &=& \displaystyle\int\limits_{\Omega} \psi_k \frac{\partial\phi_i}{\partial x}, \\
    B_y &=& [b_{y,ki}], & b_{y,kj} &=& \displaystyle\int\limits_{\Omega} \psi_k \frac{\partial\phi_j}{\partial y},
  \end{array}
\end{equation}
where $\{(\phi_1,0)^T, \dots, (\phi_m,0)^T, (0,\phi_1)^T, \dots, (0, \phi_m)^T\}$ is the basis of $\mathbf{X}_0^h$ and $\{\psi_i\}$ is the
basis of $M^h$. \RT{The case $d = 3$ is similar with the matrix $A_S$ being a $3
\times 3$ block diagonal matrix with matrix $\widehat{A}_S$ on the diagonal, and the matrix
$B$ having an additional directional component $B_z$ computed similarly to
$B_x$ and $B_y$.} For both $d=2$ and $d=3$ there is no coupling between velocity
components for different directions.

\begin{figure}[t]
\centering
  \includegraphics[scale=0.45]{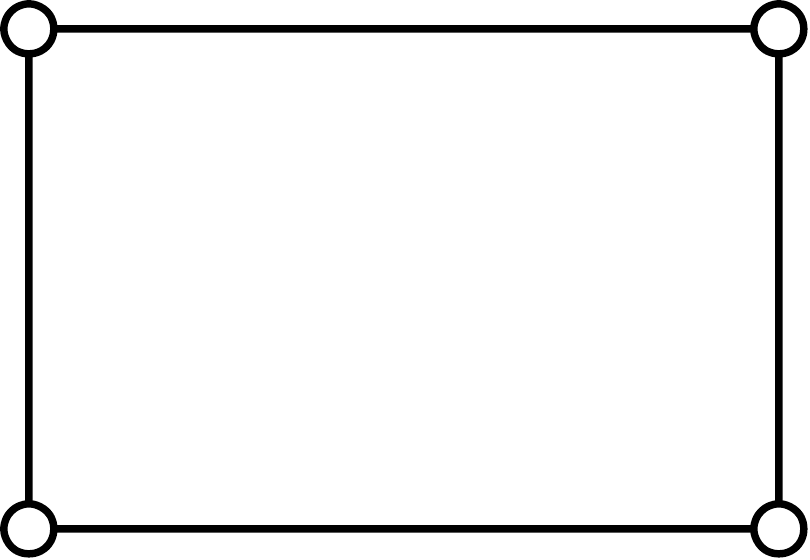} \qquad\qquad
  \includegraphics[scale=0.45]{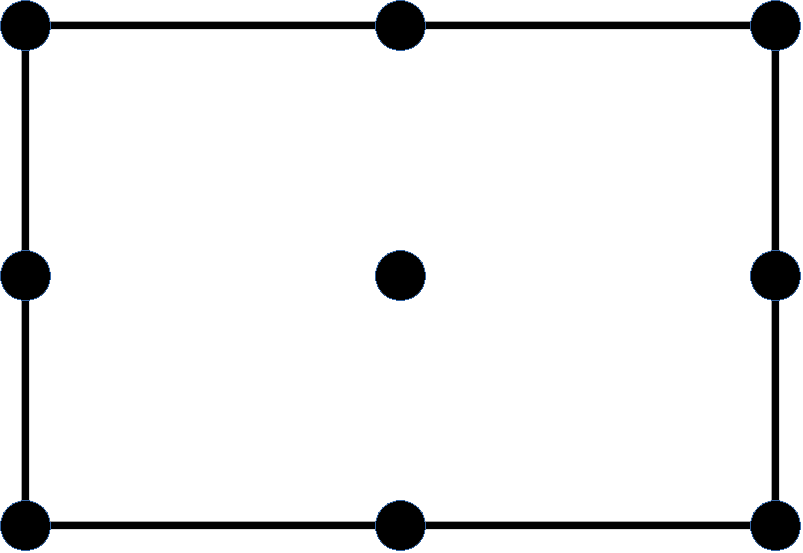}
\caption{\qq pressure (left) and velocity (right) elements.}
\label{f:q2q1_element}
\end{figure}

Using a biquadratic approximation for velocity and a bilinear approximation for
pressure in two dimensions produces a stable mixed method
called the \qq, Taylor-Hood~\cite{TaylorHood}, approximation. The
distribution of the dofs in the element is shown in Figure~\ref{f:q2q1_element}.
One important advantage of the \qq discretization is that it satisfies
\textit{inf-sup} or \textit{LBB stability}
conditions, %\improve{AP}{Cite~\cite{Elman2005}?}
and thus yields
a stable discretization~\cite{brezzi91}.  Approximations that do not satisfy this
condition require the presence of an additional stabilization term, which
introduces a stabilization matrix that replaces the zero block
in~\eqref{eq:stokes_block_system}.
%However, the \qq discretization does not require such a stabilization matrix.

\subsection{Discretization of the Navier-Stokes equations}
The steady-state Navier-Stokes equation system
\begin{equation}
  \begin{array}{rcl}
    -\nu \nabla^2 \bm{u} + \bm{u} \cdot \nabla \bm{u} + \nabla p &=& \bm{f} \quad\mbox{in}\;\Omega, \\
    \nabla \cdot \bm{u} &=& 0 \quad\mbox{in}\;\Omega,
  \end{array}
\label{eq:ns_system}
\end{equation}
with $\nu$ being a kinematic viscosity coefficient is complemented with boundary conditions given by
\begin{equation}
  \bm{u} = \bm{w} \;\mbox{on}\; \partial\Omega_D, \qquad
  \nu\frac{\partial\bm{u}}{\partial\bm{n}} - \bm{n}p = 0
  \;\mbox{on}\; \partial\Omega_N.
\label{eq:nstokes_bc}
\end{equation}
%\unsure{AP}{Why is it 0 instead of $s$ in Neumann b.c. in Elman's book?}
Linearization and discretization of this system using a Picard iteration produces the Oseen system
\begin{equation}
  \begin{pmatrix} A_{NS} & B^T \\ B & O \end{pmatrix}
  \begin{pmatrix} \mathbf{u} \\ \mathbf{p} \end{pmatrix} =
  \begin{pmatrix} \mathbf{f}_{\mathbf{u}} \\ \mathbf{f}_{\mathbf{p}} \end{pmatrix},
\label{eq:ns_block_system}
\end{equation}
with $A_{NS} = \nu A_S + K$, where $K$ is a block diagonal matrix with matrix
$\widehat{K}$ on the diagonal,
\begin{equation}
  \widehat{K} = [k_{ij}], \qquad k_{ij} = \displaystyle\int\limits_{\Omega} (\mathbf{u}_h
  \cdot \nabla \phi_j) \cdot \phi_i.
\end{equation}
Here, $\mathbf{u}_h$ is the current estimate of the discrete velocity and
once again a biquadratic basis is used for velocity while a bilinear
basis is used for pressures to arrive at a \qq discretization of the
Navier-Stokes equations.
A quantitative measure of the relative contributions of the convection and
diffusion terms is defined by the \textit{Reynolds number}, which can be written
as
\begin{equation}
  \mathcal{R} = \frac{UL}{\nu},
\end{equation}
where $U$ ($L$) is the characteristic velocity (length).

\section{Coarse level stabilization}\label{s:example}
Discretization stability gives some measure of how well-posed the discrete
problem is in comparison with the original continuous problem.
It is well known in numerical analysis that the overall quality of a
discretization can be expressed in terms of both consistency and stability.
Stability concerns arise in the multigrid context due to the construction
of the discretization operators associated with different resolutions
within the hierarchy.
Unfortunately, the stability of the finest level
discretization does not guarantee that the coarse level
operators of a multigrid hierarchy will also be stable.
Instability of any coarse level operator typically leads to dramatic
degradation in the overall multigrid convergence rates. This is often
due to the fact that standard multigrid smoothers frequently diverge
when applied to an unstable operator. Two well known areas where
discretization stability concerns arise include PDEs with highly
convective terms (even scalar PDEs) and saddle point
systems.  One possible remedy for highly convective systems lies in the use
of Petrov-Galerkin style projections~\cite{BROOKS1982199}. In this
paper we mimic this idea %in the AMG context
by employing two separate
EMIN-AMG invocations to generate interpolation and to generate
restriction~\cite{Olson2011,WiTuWaGe2012}.  To address saddle point
stability concerns, which is the focus of this paper, we seek to encourage
a relationship between the grid transfers for velocities and for pressures
that mimics that of the \qq discretization.

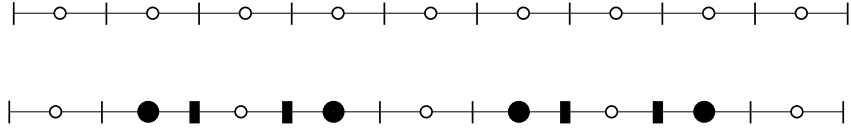
\begin{figure}[t]
\centering
\ifray
  \includegraphics[scale=.5]{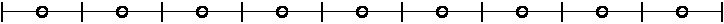}
\vskip .4in
  \includegraphics[scale=.5]{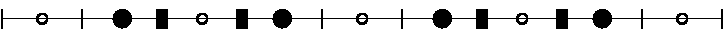}
\else
  \scalebox{.8}{\begin{tikzpicture}[y=-1cm,scale=0.6]

% objects at depth 50:
\draw[thick,black,fill=white] (4.445,10.16) circle (0.15875cm);
\draw[thick,black,fill=white] (6.985,10.16) circle (0.15875cm);
\draw[thick,black,fill=white] (9.525,10.16) circle (0.15875cm);
\draw[thick,black,fill=white] (12.065,10.16) circle (0.15875cm);
\draw[thick,black,fill=white] (14.605,10.16) circle (0.15875cm);
\draw[thick,black,fill=white] (17.145,10.16) circle (0.15875cm);
\draw[thick,black,fill=white] (19.685,10.16) circle (0.15875cm);
\draw[thick,black,fill=white] (22.225,10.16) circle (0.15875cm);
\draw[thick,black,fill=white] (24.765,10.16) circle (0.15875cm);
\draw[thick,black] (3.175,10.4775) -- (3.175,9.8579);
\draw[thick,black] (5.715,10.4775) -- (5.715,9.8679);
\draw[thick,black] (8.255,10.4775) -- (8.255,9.8679);
\draw[thick,black] (10.795,10.4775) -- (10.795,9.8679);
\draw[thick,black] (13.335,10.4775) -- (13.335,9.8679);
\draw[thick,black] (15.875,10.4775) -- (15.875,9.8679);
\draw[thick,black] (18.415,10.4775) -- (18.415,9.8679);
\draw[thick,black] (20.955,10.4775) -- (20.955,9.8679);
\draw[thick,black] (23.495,10.4775) -- (23.495,9.8679);
\draw[thick,black] (26.035,10.4775) -- (26.035,9.8679);

% select the background layer
\begin{pgfonlayer}{bg}
  \draw[black] (3.175,10.16) -- (26.035,10.16);
\end{pgfonlayer}

\end{tikzpicture}%
}
\vskip .35in
  \scalebox{.8}{\begin{tikzpicture}[y=-1cm,scale=0.6]

% objects at depth 50:
\draw[thick,black,fill=white] (4.445,10.16) circle (0.15875cm);
\draw[thick,black,fill=white] (6.985,10.16) circle (0.15875cm);
\draw[thick,black,fill=white] (9.525,10.16) circle (0.15875cm);
\draw[thick,black,fill=white] (12.065,10.16) circle (0.15875cm);
\draw[thick,black,fill=white] (14.605,10.16) circle (0.15875cm);
\draw[thick,black,fill=white] (17.145,10.16) circle (0.15875cm);
\draw[thick,black,fill=white] (19.685,10.16) circle (0.15875cm);
\draw[thick,black,fill=white] (22.225,10.16) circle (0.15875cm);
\draw[thick,black,fill=white] (24.765,10.16) circle (0.15875cm);
\draw[line width=4.8bp,black] (6.985,10.16) circle (0.15875cm);
\draw[line width=4.8bp,black] (12.065,10.16) circle (0.15875cm);
\draw[line width=4.8bp,black] (17.145,10.16) circle (0.15875cm);
\draw[line width=4.8bp,black] (22.225,10.16) circle (0.15875cm);
\draw[thick,black] (3.175,10.4775) -- (3.175,9.8552);
\draw[thick,black] (5.715,10.4775) -- (5.715,9.8679);
\draw[thick,black] (13.335,10.4775) -- (13.335,9.8679);
\draw[thick,black] (15.875,10.4775) -- (15.875,9.8679);
\draw[thick,black] (23.495,10.4775) -- (23.495,9.8679);
\draw[thick,black] (26.035,10.4775) -- (26.035,9.8679);
\draw[line width=4.8bp,black] (8.255,10.4775) -- (8.255,9.8679);
\draw[line width=4.8bp,black] (10.795,10.4775) -- (10.795,9.8679);
\draw[line width=4.8bp,black] (18.415,10.4775) -- (18.415,9.8679);
\draw[line width=4.8bp,black] (20.955,10.4775) -- (20.955,9.8679);

% select the background layer
\begin{pgfonlayer}{bg}
  \draw[black] (3.175,10.16) -- (26.035,10.16);
\end{pgfonlayer}

\end{tikzpicture}%
}
\fi
\caption{One dimensional mesh.
Top: fine pressure (velocity) unknowns denoted by circles (vertical lines).
Bottom: coarse pressure (velocity) unknowns denoted by filled circles (boxes).}
\label{f:macmesh}
\end{figure}

\begin{figure}[t]
\centering
\ifray
  \includegraphics[scale=.2]{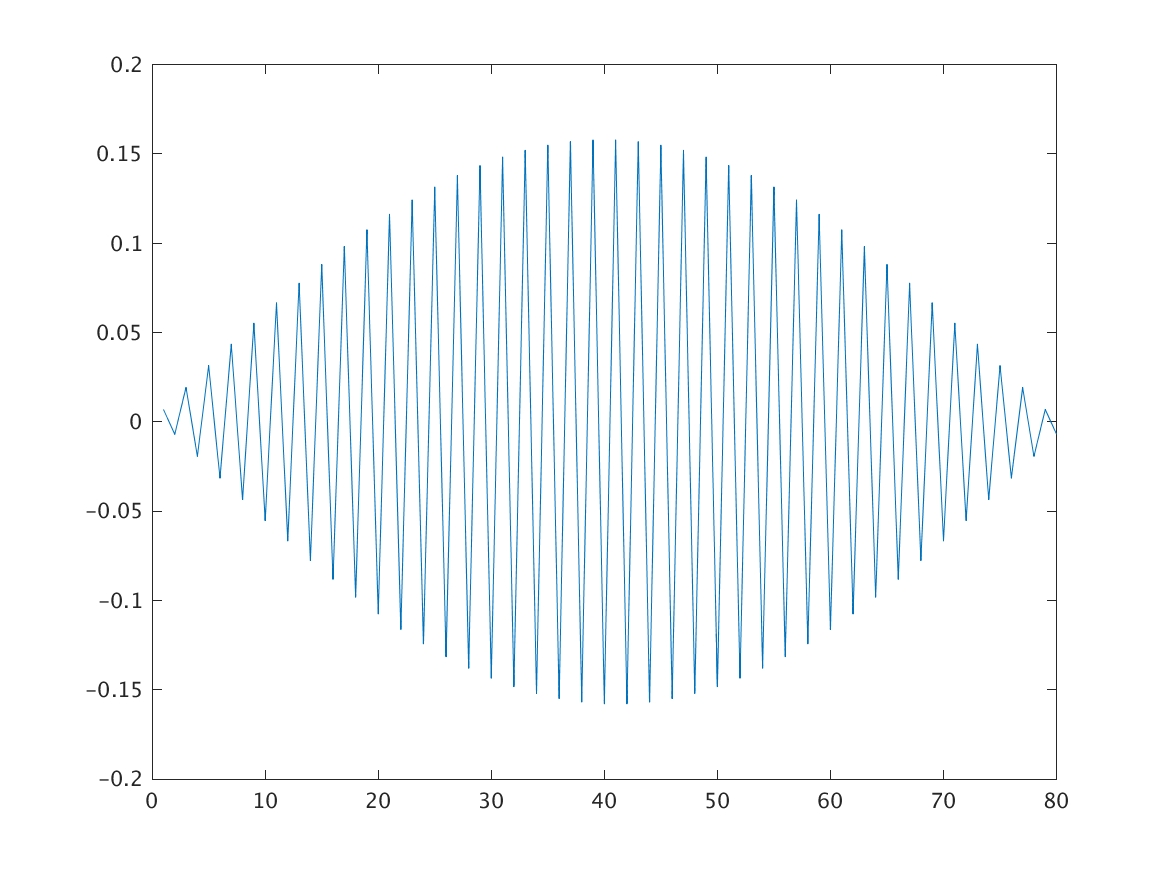}
\else
  \scalebox{0.9}{% This file was created by matlab2tikz.
\begin{tikzpicture}[scale=0.6]

\begin{axis}[%
width=5.291in,
height=4.173in,
at={(0.887in,0.563in)},
scale only axis,
separate axis lines,
every outer x axis line/.append style={black},
every x tick label/.append style={font=\color{black}},
x tick label style={
  font=\Large,
},
xmin=0,
xmax=80,
every outer y axis line/.append style={black},
every y tick label/.append style={font=\color{black}},
y tick label style={
  /pgf/number format/.cd,
  fixed,
  fixed zerofill,
  precision=2,
  /tikz/.cd,
  font=\Large,
},
ymin=-0.2,
ymax=0.2,
axis background/.style={fill=white}
]
\addplot [color=blue,solid,forget plot]
  table[row sep=crcr]{%
1	0\\
2	0.0124054722828585\\
3	-0.0124054722828583\\
4	0.0247344607038551\\
5	-0.024734460703855\\
6	0.0369109529495693\\
7	-0.0369109529495693\\
8	0.0488598768962119\\
9	-0.0488598768962118\\
10	0.0605075634542314\\
11	-0.0605075634542315\\
12	0.0717822007627489\\
13	-0.071782200762749\\
14	0.0826142769335597\\
15	-0.0826142769335604\\
16	0.0929370086150406\\
17	-0.0929370086150417\\
18	0.102686752733719\\
19	-0.10268675273372\\
20	0.111803398874983\\
21	-0.111803398874985\\
22	0.12023073988378\\
23	-0.120230739883782\\
24	0.127916818400416\\
25	-0.127916818400419\\
26	0.134814247194955\\
27	-0.134814247194958\\
28	0.14088050132525\\
29	-0.140880501325254\\
30	0.146078180317356\\
31	-0.14607818031736\\
32	0.15037523875188\\
33	-0.150375238751884\\
34	0.153745183834659\\
35	-0.153745183834663\\
36	0.156167238733635\\
37	-0.156167238733639\\
38	0.157626470674942\\
39	-0.157626470674946\\
40	0.158113883008416\\
41	-0.158113883008421\\
42	0.157626470674943\\
43	-0.157626470674949\\
44	0.156167238733638\\
45	-0.156167238733644\\
46	0.153745183834663\\
47	-0.15374518383467\\
48	0.150375238751887\\
49	-0.150375238751893\\
50	0.146078180317363\\
51	-0.146078180317369\\
52	0.140880501325258\\
53	-0.140880501325265\\
54	0.134814247194963\\
55	-0.13481424719497\\
56	0.127916818400424\\
57	-0.127916818400432\\
58	0.120230739883788\\
59	-0.120230739883796\\
60	0.111803398874991\\
61	-0.111803398874999\\
62	0.102686752733726\\
63	-0.102686752733735\\
64	0.0929370086150478\\
65	-0.0929370086150564\\
66	0.0826142769335654\\
67	-0.0826142769335743\\
68	0.0717822007627528\\
69	-0.0717822007627619\\
70	0.060507563454234\\
71	-0.0605075634542433\\
72	0.0488598768962123\\
73	-0.0488598768962218\\
74	0.0369109529495685\\
75	-0.0369109529495781\\
76	0.0247344607038528\\
77	-0.0247344607038625\\
78	0.0124054722828548\\
79	-0.0124054722828645\\
80	-4.83234796483599e-15\\
};
\end{axis}
\end{tikzpicture}%
}
\fi
\caption{Eigenvector associated with second smallest eigenvalue of
Schur complement of projected system with coarse points as in Figure~\ref{f:macmesh}.}
\label{f:oscillatory}
\end{figure}
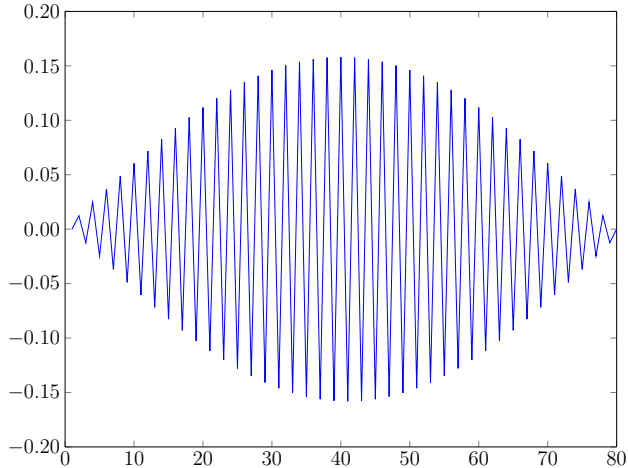

The paper \cite{wabro2004} discusses cases where instability
of AMG coarse grid operators can appear and how this instability can
lead to significant multigrid convergence degradation. Here, we give
a rather elementary example to illustrate how easily instability can arise
when the coarsening strategy for pressures and velocities is
somewhat inconsistent. Specifically, consider the following simplified matrix
which can be viewed as a marker-and-cell style approach~\cite{HarlowWelch}
to a one dimensional PDE:
\begin{equation} \label{mac system}
\begin{pmatrix}
~I & B^T \\
-B & O   \\
\end{pmatrix}   \begin{pmatrix} u \\ p \end{pmatrix} =
\begin{pmatrix} f \\ 0 \end{pmatrix},
\end{equation}
where
\begin{equation} \label{mac submatrix}
B^T = \begin{pmatrix}
-1 &  1 &                         \\
   &  . & . &                     \\
   &    & . & .  &                \\
   &    &   & .  & .  &           \\
   &    &   &    & -1 &  1        \\
   &    &   &    &    & -1 & 1
\end{pmatrix}.
\end{equation}
$ B^T \in \mathbb{R}^{n \times n\hskip-.01in-\hskip-.01in1}$ can be viewed as a
discrete gradient operator and $B$ as a divergence operator where pressure
(velocity) unknowns are defined at circles (vertical lines\footnote{The two
  velocity end points correspond to Dirichlet boundary conditions that have been
eliminated from \eqref{mac system}.}) in the top image of
Figure~\ref{f:macmesh}. Here, $n$ denotes the number of fine level pressure
dofs.
The Schur complement matrix, $B B^T$, is identical to that obtained
by a standard three point central difference discretization of the
Poisson operator with Neumann boundary conditions.
Now consider a coarse mesh (bottom image of Figure~\ref{f:macmesh})
in conjunction with linear
interpolation\footnote{
Fine points between two coarse points use linear interpolation.
%Fine points aligned with coarse points use injection.
Fine pressures at end points injection from the closest coarse pressure (preserving
constant null space). End point fine velocities use linear interpolation with
Dirichlet end point (assuming value of zero at Dirichlet point).}
operators: $P_v$ and $P_p$.
Thus, $P_p$ is a fairly standard geometric interpolation operator, and
$P_v$ is close to standard geometric interpolation with the one
exception that the choice of coarse points is somewhat sub-optimal.
In this case, the stencil of the projected gradient, $P_v^T B^T P_p$,
resembles a centered difference scheme, and the resulting projected
operator is
\begin{equation} \label{projected mac}
\begin{pmatrix}
  ~P_v^T P_v   & P_v^T B^T P_p \\
  -P_p^T B P_v & O             \\
\end{pmatrix}.
\end{equation}
The interior stencil of the Schur
complement of the projected operator, $S = P_p^T B P_v (P_v^T P_v)^{-1} P_v^T B^T P_p$, only weakly couples even and odd coarse grid points. In fact, the
related operator $ \widehat{S} = P_p^T B P_v P_v^T B^T P_p$ has completely
decoupled even and odd points,
%resembles a one dimensional
%interior five point Laplace stencil such that odd and even points are no longer
%coupled,
i.e. $ \widehat{S}_{i,i\pm2} = -1,  \widehat{S}_{i,i\pm1} = 0,  \widehat{S}_{i,i} = 2$.
Both $S$ and $\widehat{S}$ are unstable Laplacians as indicated by the highly
oscillatory eigenvector associated with the second smallest eigenvalue
of $S$ shown in Figure~\ref{f:oscillatory}.
An oscillatory mode associated with a small eigenvalue of $S$ (or equivalently
a large eigenvalue of $S^{-1}$) is an indication
that we can generally expect
$S^{-1} v$
%\fix{AP}{$S^{-1}$ cannot be applied to $f$ due to different dimensions. Should the rhs in~\eqref{mac system} be (0, f) instead?}
to significantly amplify oscillatory components of a vector $v$.
%to be highly oscillatory.
This instability is purely a function of the poor choice for the
velocity coarse grid points. In fact, no instability is present
in the resulting Schur complement if instead the coarse velocities
are associated with mid-points between the pressures and linear
interpolation is used.

$S$'s instability is likely to be problematic if AMG is applied to $S$
to further coarsen pressures and when a relaxation scheme is
applied to \eqref{projected mac}.  Further,
we note that many specialized relaxation schemes developed for PDE systems
rely on sub-matrix properties mirroring
those of the corresponding PDE operators
(e.g., $P_v^T B^T P_p$ mirroring a gradient operator).
When this no longer holds,
the specialized relaxation schemes may not smooth errors as expected.

\section{AMG Grid Transfers for \qq discretizations}\label{s:algorithm}
This section proposes an algorithm for generating grid transfers
while the next section describes appropriate relaxation procedures.
Together, these two components fully define an AMG method
when used in conjunction with a Petrov-Galerkin projection process.
\RT{The specific technique leverages the fact that a subset of the 
velocity unknowns are co-located with the pressure unknowns
for \qq discretizations, which would not be the case for all types of mixed
finite element discretizations.}
Let $A$ be a given full system matrix, $A = A_{NS}$ for concreteness.
Algorithm~\ref{a:AMG_setup} summarizes the setup of an AMG
hierarchy given procedures for generating grid transfers and
relaxation operators.
\begin{algorithm}[t]
\caption{\func{setup\_hierarchy}$(A_\ell, x_\ell, \ell, \ell_{max})$}
\label{a:AMG_setup}
  \begin{algorithmic}
    \State $S_\ell = \mbox{\func{form\_relaxation}}(A_\ell)$
    \If {$\;\ell < \ell_{max}\;$}
      \State $P_\ell = \mbox{\func{form\_interpolation}}(A_\ell, x_\ell)$
      \State $R_\ell = \mbox{\func{form\_restriction}}(A_\ell, x_\ell, P_\ell)$
       \State $A_{\ell+1} = R_\ell A_\ell P_\ell $
       \State $x_{\ell+1} = \mbox{\func{project\_coordinates}}(x_\ell, R_\ell)$
       \State  \func{setup\_hierarchy}$(A_{\ell+1}, x_{\ell+1}, \ell+1, \ell_{max})$
    \EndIf
  \end{algorithmic}
\end{algorithm}
In Algorithm~\ref{a:AMG_setup}, $\ell_{max}$ is the number of levels in a
multigrid hierarchy, $\ell = 1$ corresponds to the finest level, the matrix
$P_\ell$ defines interpolation from hierarchy level $\ell+1$ to level $\ell$
while $R_\ell$ defines restriction from $\ell$ to $\ell+1$. $A_\ell$ is the
``discretization''
%\unsure{AP}{Is there a better choice of words?}
matrix on the $\ell$th level (with $A_0 = A$ being the finest level matrix), and
$S_\ell$ is the smoother that is used within a relaxation procedure of the form
$$
u_\ell^{k+1} =  u_\ell^k + S_\ell ( b_\ell - A_\ell u_\ell^k ),
$$
where $u_\ell^k$ is an approximation to the solution on level $\ell$, and $b_\ell$
is the right hand side on the same level.
As our algorithms rely on coordinate
locations $x_\ell$ (but not other mesh information), this is
illustrated in the code fragment.

The discussion focuses on construction of interpolation matrices as the
procedure for generating restriction operators only includes some modest
differences. We also center the exposition on the principal themes while
omitting some detailed graph heuristics. % in the Appendix.
These details
generally involve tie-breaking and weighing tradeoffs
between including additional coarse points that might improve convergence,
but increase iteration cost (and memory) due to increased coarse level fill-in.
Many AMG codes have similar types of heuristics, and these
are not the focus of this article.
For the incompressible
Navier-Stokes equations, our prolongators have the block diagonal form
$$
    P_\ell = \begin{pmatrix} P^{(v)}_\ell & O \\ O & P^{(p)}_\ell \end{pmatrix},
$$
which is the same as \eqref{eq:general_block_p} with slightly different notation
to simplify the presentation. The basic idea is
to apply an AMG algorithm to first generate a pressure prolongator $P^{(p)}_\ell$ followed by a
procedure that generates a velocity prolongator $P^{(v)}_\ell$ using coarse points that are
\textit{consistent} with those used for $P^{(p)}_\ell$. By consistent, we
mean that any coarse pressure point has co-located coarse velocity points
as well as some notion of coarse velocity mid-points in between the
co-located coarse points.
The initial focus in our discussion is on coarse point selection and
the choice of grid transfer sparsity patterns as these are the
primary inputs to an energy minimization AMG algorithm. In addition
to a consistent coarse point selection, we also aim to produce a
multigrid algorithm with relatively sparse coarse grid operators.
That is, we want coarsening rates and sparsity patterns
to be defined in a way such that the matrix projection
$ R_\ell A_\ell P_\ell $ is not significantly denser than $A_\ell $.
This requires some special considerations when biquadratic basis functions
are employed on the finest level.

Algorithm~\ref{a:build_block_P} illustrates a general strategy.
To simplify notation, level sub-scripts are omitted.
% from the discussion.
\begin{algorithm}[t]
\caption{$P = \mbox{\func{form\_interpolation}}(A, x)$}
\label{a:build_block_P}
  \begin{algorithmic}
    \State $ \left[\widetilde{A}^{(v)}, \widetilde{A}^{(p)}\right] = \mbox{\func{form\_aux\_block\_diagonal}}(A)$\\[    ]
    \State $ C^{(p)} = \mbox{\func{find\_coarse\_pressure}}(\widetilde{A}^{(p)})$
    \State $ N^{(p)} = \mbox{\func{form\_pressure\_pattern}}(\widetilde{A}^{(p)}, C^{(p)})$
    \State $ P^{(p)} = \mbox{\func{form\_pressure\_prolongator}}(\widetilde{A}^{(p)}, C^{(p)}, N^{(p)})$\\[    ]
    \State $ C^{(v)} = \mbox{\func{form\_coarse\_velocities}}(\widetilde{A}^{(p)}, C^{(p)}, N^{(p)}, x)$
    \State $ N^{(v)} = \mbox{\func{form\_velocity\_pattern}}(\widetilde{A}^{(v)}, C^{(v)})$
    \State $ P^{(v)} = \mbox{\func{form\_velocity\_prolongator}}(\widetilde{A}^{(v)}, C^{(v)}, N^{(v)})$
  \end{algorithmic}
\end{algorithm}
% \info[inline]{AP}{To be completely fair, in Algorithm~\ref{a:build_block_P} it should be
% $ C^{(p)} = \mbox{\func{form\_coarse\_pressures}}(\widetilde{A}^{(p)}, x)$
% and
% $ C^{(v)} = \mbox{\func{form\_coarse\_velocities}}(\widetilde{A}^{(v)}, C^{(p)},
% N^{(p)}, x)$. Specifically, when we convert some $F$ vertices into $C$
% vertices, we use graph information, specifically, the velocity graph
% information. But talking about this would complicate the paper.}
First, submatrices $\widetilde{A}^{(v)}$ and $\widetilde{A}^{(p)}$
are formed corresponding to velocity-velocity and presure-pressure coupling;
subsequent AMG procedures will be applied to these matrices. While $A_{NS}$
could be used for $\widetilde{A}^{(v)}$, it may be advantageous to
filter out {\it weak} couplings~\cite{BrMcRu1984,Ruge1987}. Obviously,
more care is required in defining $\widetilde{A}^{(p)}$ as the $(2,2)$ block
of the PDE system matrix is identically zero. Coarse points are then defined.
Specifically, $C^{(p)}$ and $C^{(v)}$ are sets of indices corresponding to
the fine level degrees-of-freedom (pressures and velocities respectively)
that will be represented on the next coarsest level.  Formally,
$C^{(p)} \subset \mathbb{Z}_{n^{(p)}}$ and
$C^{(v)} \subset \mathbb{Z}_{n^{(v)}}$,
where $n^{(p)}$ ($n^{(v)}$) is  the number of fine grid pressures (velocities)
and $\mathbb{Z}_k$ denotes the set of integers between $1$ and $k$ inclusive.
Sparsity pattern matrices for the grid transfers are then constructed.
For example, $N^{(p)}$ is a binary matrix indicating that
$P^{(p)}_{ij}$ is nonzero
only if $N^{(p)}_{ij}$ is one. Thus, the $j$th column of $N^{(p)}$
defines the support of the $j$th interpolating basis function.
Finally, actual interpolation weights must be calculated. Given that
coarse points and sparsity patterns will be chosen in somewhat nonstandard
ways to preserve certain discretization features, it is important
that a flexible algorithm is used to define interpolation weights.
That is, the algorithm must not rely on restrictive assumptions concerning
the distances between coarse points or the grid transfer sparsity patterns.
For this, we employ an energy-minimizing framework (EMIN-AMG)~\cite{Olson2011},
which is described further below.

%The main idea with energy minimization AMG is that ...
%
%\begin{verbatim}
%[ Perhaps we can limit the null space description and somehow just say
%  in our case N_C is N are both binary (0/1 matrices) that specify
%that coarse constant vectors interpolate to fine constant vectors.
%\end{verbatim}

\subsection{Coarse points and sparsity patterns for pressure grid
transfers}\label{s:ptransfer}

%Coarsening and sparsity pattern algorithms are applied to
The $\widetilde{A}^{(p)}$ matrix
%and so this
must first be defined.  A natural choice would be to
take $\widetilde{A}^{(p)} = B B^T$ as this leads to a pressure Poisson operator,
which often plays a central role within incompressible Navier-Stokes
calculations. One complication, however, is that the stencil for this
pressure Poisson operator is somewhat non-standard. In particular,
the use of biquadratic basis functions for the velocities implies that
stencils associated with $A_{NS}, B $ and consequently $B B^T$ are
wider/fuller than those obtained with standard first order finite differences,
finite volumes, or finite elements.
First order schemes typically give rise
to fairly compact stencils involving only mesh points that immediately
surround a central point.
For example, the \qq interior stencils
associated with $B B^T$ are 25 point stencils (as opposed to standard
nine point stencils) when the mesh is a regular two-dimensional structured
grid. To encourage stencil widths that more closely resemble standard
first order discretizations, $\widetilde{A}^{(p)}$ is defined by dropping small
entries in $Z = B B^T$ and lumping dropped entries to the diagonal so the sum of
matrix entries within a row is preserved. This dropping essentially
corresponds to removing nonzero entries that satisfy
\begin{equation}
|z_{ij}| \le \tau_1  \sqrt{| z_{ii} z_{jj} |},
\label{e:filtering}
\end{equation}
with $\tau_1$ being a user provided parameter.
%See Appendix~\ref{s:appendix_dropping} for details.

%As depicted in Algorithm~\ref{a:build_block_P}, $P^{(p)}$ does not
%depend on the velocities. Thus, it might be possible to apply
%a standard AMG algorithm to define $P^{(p)}$. However, standard AMG
%algorithms are best suited to first order
%discretization schemes.

%Due to potential higher order discretization issues,
%we choose to employ energy minimization AMG to generate
%$P^{(p)}$ applying it to operators that more closely resemble those
%produced by standard first order discretizations.
% This gives rise to nine point
% stencils if used with a threshold of $10^{-?}$ and applied to $B B^T$
% when the underlying mesh is a two-dimensional structured grid.

The other non-standard aspect of our AMG procedure for $P^{(p)}$
is that we orient the algorithm so that coarse pressure points are
more distant than normal. In particular, classical AMG targets
a set of coarse points that are distance two from each other in
the associated matrix graph. Smoothed aggregation targets a set of
aggregate root points (the aggregation counterpart to coarse points)
that are distance three from each other.  In our algorithm, we
target distance four coarse points in the graph associated with
$\widetilde{A}^{(p)}$.  This choice is again driven by stencil widths
within high order discretizations and avoiding excessive fill-in
during the Petrov-Galerkin projection, $R_\ell A_\ell P_\ell$.
When distance four coarse points are used in
conjunction with the sparsity pattern to be discussed, the projection of
$B$ (and $A_{NS}$) to coarser levels produce little additional fill-in
(i.e. the average number of nonzeros per row does not rise appreciably).
\RT{However, more traditional coarsening rate of three might be worth
further exploration given the low complexity rates demonstrated in our
numerical results with the coarsening rate of four.}
%This is discussed further below.

\begin{algorithm}[t]
\caption{$C^{(p)} = \mbox{\func{find\_coarse\_pressures}}(A, x)$}
\label{a:PressureStep}
  \begin{algorithmic}
    \State $[V, E] = \mbox{\func{form\_graph}}(A)$ \Comment{ $A \in
    \mathbb{R}^{n \times n}$, $V$ is the set of vertices,  $E$ is the set of edges}
    \State $C^{(p)} = \varnothing ~;~ F^{(p)} = \varnothing $ \Comment{$C^{(p)} (F^{(p)})$ is a set of coarse (fine) vertices}
    \State $S = \{S_i\}_i, \;S_i = \varnothing \hskip .2in \forall i \in \mathbb{Z}_n$
    \State $Cand = \varnothing$  %\Comment{Arbitrary initial point}
    \State $k = \mbox{arbitrary vertex}$
    \While{ $F^{(p)} \cup C^{(p)} \subsetneq \mathbb{Z}_n$ }
      \State $D_3 = \{ ~j ~|~ dist(V_k,V_j, E) \le 3~ \}$ \
      \State $D_4 = \{ ~j ~|~ dist(V_k,V_j, E) = 4~ \}$
      \State $S_j = S_j \cup k \hskip .2in \forall j \in D_3 $ \Comment{\RT{$S_j$ updated to include all coarse points within distance 3 of $V_j$}}
      \State $D_3 = D_3 \setminus k$
      \State $C^{(p)} = C^{(p)} \cup k $
      \State $F^{(p)} = F^{(p)} \cup D_3$
%      \State $Cand = Cand \setminus D_3$%:= \{ ~j \in Cand~ | ~j \notin D_3 \}$
      \State $Cand = Cand \cup D_4 $
      \State $Cand = Cand \setminus (F^{(p)} \cup C^{(p)}) $
      \State $[h_1, h_2] = \mbox{\func{update\_heuristics}}(V, E, S, x, k, h_1, h_2)$
      \State $k = \argmin_{j \in Cand} h_1(j)$ \Comment{Best $Cand$ or arbitrary unmarked}
      \State \Comment{vertex if $Cand = \varnothing$~~~~~~~~~~~~~~~~~~~~~~}
    \EndWhile
    \State $[C^{(p)}, F^{(p)}, S] = \mbox{\func{find\_extra\_dist3\_Cpoints}}(V, E, S, C^{(p)}, F^{(p)}, x, h_1, h_2)$
%    \State $T = \{ j ~|~ S_j \le 2  \mbox{ and } j \in F \}$
%    \While{ $T \neq \varnothing$ }
%        \State $k = $\func{chooseElement}$(T)$
%    \EndWhile
  \end{algorithmic}
\end{algorithm}

The general coarse point selection
algorithm is given in Algorithm~\ref{a:PressureStep}.
The idea is to classify all $ n^{(p)} $ vertices as either coarse points
(represented on next coarse grid) or fine points (not represented on next coarse grid).
In the AMG literature, these sets are respectively referred to as the
$C$- and $F$-point sets, and the classification as a $C/F$-splitting. In
this paper, we use $C^{(p)}$ and $F^{(p)}$ ($C^{(v)}$ and $F^{(v)}$) to refer to
these two sets for pressure (velocity) vertices. Let $G = (V, E)$ be a graph of
matrix $\widetilde{A}^{(p)}$.
The classification is performed by first selecting an initial vertex, $k \in V$,
and including it in $C^{(p)}$. All vertices within distance three from $k$,
$dist(V_k, V_j,E) \le 3$,
(excluding $k$) are added to $F^{(p)}$.
Here, distance refers to the minimum number of edge traversals
(defined by $E$) required to travel from one vertex to another.
A candidate set, $Cand$,
is also introduced to encourage future coarse points to be distance four from
existing $C^{(p)}$. Specifically, $Cand$ is the set of all vertices
that are distance four from any $C^{(p)}$ and have not been already
included in $C^{(p)}$ or $F^{(p)}$.
The next $C$-point is selected
giving priority to vertices in $Cand$ using heuristics that
encourage the close {\it packing} of coarse points
(by choosing a vertex with the smallest average Euclidean distance to
existing  $C^{(p)}$ points).
The above procedure is repeated until all vertices are classified.
Upon completion, all vertices in $C^{(p)}$ are at least distance four from
each other. All vertices in $F^{(p)}$ are within distance
three from at least one vertex in $C^{(p)}$.  This may lead to situations
where some $F^{(p)}$ vertices are not well covered by $C^{(p)}$ vertices, e.g.
only within distance three from just a single $C^{(p)}$. This
type of issue occurs in most algebraic multigrid codes (including our
smoothed aggregation library) and can degrade convergence rates.
To minimize these effects, heuristics are generally used
%In our case,
%Appendix~\ref{s:appendix_extra_dist3_Cpoints} provides details of the
%\func{find\_extra\_dist3\_Cpoints} that
to convert some $F^{(p)}$ vertices to $C^{(p)}$ vertices, targeting those
not well covered by current $C^{(p)}$ vertices, i.e. those with a small size of
set $S_j$ \RT{where for every $j \in F^{(p)}$ the set $S_j$ corresponds to all $C^{(p)}$ vertices that are
within graph distance three from vertex $V_j$ ($j \in \mathbb{Z}_n$)}.
Specifically, heuristics first convert $j \in F^{(p)}$ when $|S_j| = 1$
and the Euclidean distance to the one $C^{(p)}$ point is
large.
Heuristics also examine points where $|S_j| = 2$
and converts them to $C^{(p)}$ vertices only if the average Euclidean and
the average graph distance to the $C^{(p)}$ vertices is large and if
$j$ is far from the line segment connecting the two $C^{(p)}$ vertices.
Typically, a small percentage of $C^{(p)}$ vertices
are chosen by this algorithm. %Appendix~\ref{s:appendix_heuristics} gives
Harmonic averages (and standard averages) of Euclidean distances (and graph
distances) to nearby $C^{(p)}$ vertices are stored in $h_1(j)$ (and
$h_2(j)$). These are computed by \func{update\_heuristics}.
Harmonic averages essentially favor minimum distances.
% updates
%of the Euclidean distance between a vertex $j$ and any new  $C^{(p)}$ vertices
%that encourages close {\it
%packing}.
%The function
%$h_1(j)$ so that it corresponds

To define $P^{(p)}$'s sparsity pattern the sets $S_i $ are
used.\footnote{The $S_i$ sets are also used by \func{find\_extra\_dist3\_Cpoints}.}
The sparsity pattern is then given
by
$$
N^{(p)}_{ij} = \begin{cases}
1,  & \mbox{if $j \in S_i$}, \\
0,  & \mbox{otherwise}.
\end{cases}
$$

\begin{figure}[t]
\centering
  \includegraphics[scale=.30]{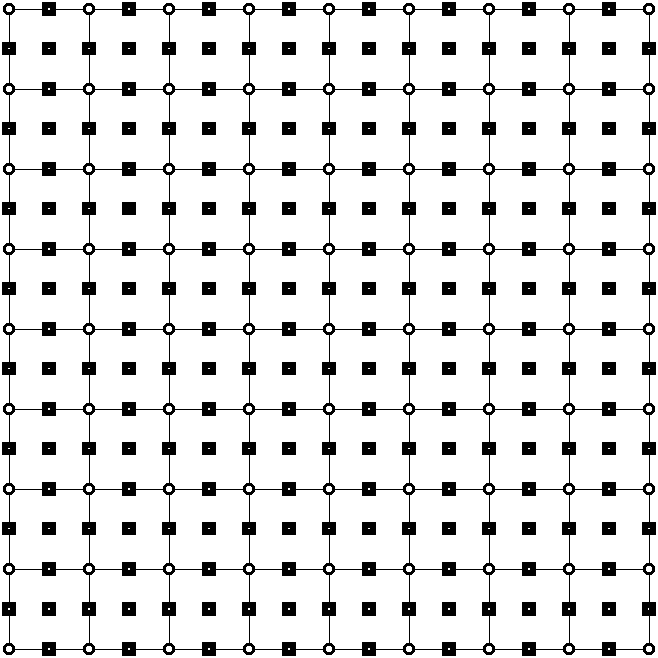}~~~~~~~~~~~
  \includegraphics[scale=.30]{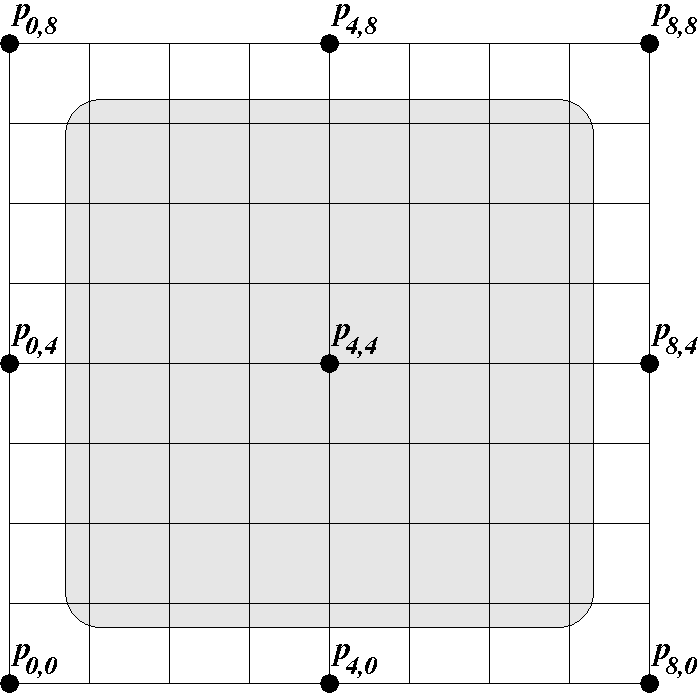}
\caption{Left: \qq discretization where circles  (boxes) are pressure and
  co-located velocity (non-co-located velocity) dofs.
%Boxes correspond to non-co-located velocity dofs.
Right: perfect pressure coarsening. Circles are coarse pressure dofs. Dofs
surrounded by gray rectangle belong to interpolation sparsity pattern
associated with middle coarse dof's column.}
\label{f:prescpts}
\end{figure}
Figure~\ref{f:prescpts} (right) illustrates a perfect
coarsening of the pressure dofs for the discretization depicted in
Figure~\ref{f:prescpts} (left). Here, all nearest neighbors to any
coarse pressure are exactly distance four away. This
graph distance is in terms of the $\widetilde{A}^{(p)}$ matrix.\footnote{The figure
displays the underlying mesh as opposed to the matrix graph. The
$\widetilde{A}^{(p)}$ matrix graph would also include diagonal edges within
each mesh box.} Such a perfect coarsening would not generally occur, though
our metrics encourage it. Consider, for example, the
situation where the $p_{0,0}$ is the initial coarse vertex. Here,
sub-script $i,j$ refers to the vertex that is $i$ ($j$) horizontal (vertical)
edges away from the lowest-leftmost vertex.
The $D_4$ and $Cand$ sets would include all vertices along the line between
$p_{0,4}$ and $p_{4,4}$ inclusive and between $p_{4,4}$ and $p_{4,0}$ inclusive.
To encourage perfect coarsening, the selection of a next $C$-point should
favor vertices {\it spatially} closest to the initial corner vertex so that
either $p_{0,4}$ or $p_{4,0}$ is chosen next. However, encouraging perfect
coarsening becomes complicated as additional coarse vertices are
chosen, and it is not necessary from a numerical convergence perspective.
Finally, the gray shaded region illustrates the interpolation sparsity
pattern associated with the $p_{4,4}$ vertex.
In particular, all vertices surrounded by the gray region are at most
distance three from $p_{4,4}$. Thus, these
degrees-of-freedom would be nonzero in the column of the matrix interpolation
operator associated with this $p_{4,4}$ coarse pressure.

\subsection{Coarse points and sparsity patterns for velocity grid
transfers}\label{s:vtransfer}
Each velocity component (in the different coordinate directions) is coarsened
in an identical fashion and uses the same grid transfer sparsity pattern.
The method for determining coarse velocity vertices is given in
Algorithm~\ref{a:midpoints}. The algorithm starts by taking pressure coarse
points, $C^{(p)}$, and adding pressure mid-points to define the set
$\overline{C}^{(p)}$.
The basic idea for mid-points is to first compute a set of target spatial
locations $X_i$. The number of unique target locations is based on the number of
unique $S_i$ sets. Recall,
%$S_i$ determines the $i$th row of $N^{(p)}$. That is,
$S_i$ gives all coarse vertices used (i.e. with nonzero interpolation weights)
when interpolating to the $i$th fine pressure vertex.  Each target location
corresponds to the geometric centroid (or barycenter) of the spatial locations
associated with $S_i$'s coarse points.  For each target location, we find the
spatially closest fine pressure vertex. The search for the closest pressure is
limited to the set of vertices $B_i$ which corresponds to all $F^{(p)}$ points
within distance three to {\it all} of $S_i$'s vertices.  Thus, the $B_i$
vertices define a small neighborhood surrounding the $i$th vertex.  The closest
pressure dof is chosen as a pressure mid-point only if it is not too close to an
already to an already chosen coarse pressure in $\overline{C}^{(p)}$. Closeness is
defined in a relative sense with respect to a box containing all $B_i$
and a tolerance, $tol$, shown in Algorithm~\ref{a:midpoints}.
Further, the order in which our implementation computes target mid-points
(the {\bf for} loop in  Algorithm~\ref{a:midpoints}) is such that
vertices with larger $|S_i|$ are chosen before those with smaller $|S_i|$
($|S_i|$ denotes the cardinality of the set $S_i$).
In the final step, the coarse velocity points are taken to be velocities
co-located with the chosen $\overline{C}^{(p)}$ pressures.
% Finally, co-located velocity dofs are determined for each $\overline{C}^{(p)}$ vertex.
% This mapping between pressure and velocity dofs must be provided to the algorithm and
% it must be projected on coarser levels within the algorithm.
\begin{algorithm}[t]
\caption{$C^{(v)} = \mbox{\func{find\_velocity\_Cpoints}}(S, C^{(p)})$}
\label{a:midpoints}
  \begin{algorithmic}
    \State $ \overline{C}^{(p)} =  C^{(p)}$
    \For{ $i \in F^{(p)}$ }
      \State $X_i = \left(\sum_{k \in S_i} x^{(p)}_k\right) /|S_i|$
      \Comment{Barycenter of surrounding $C$-points}
      \State $B_i = \{j ~|~ j \in F^{(p)} \mbox{ and } S_i \subseteq S_j \}$ \Comment{All $F$-points that also interpolate from $S_i$}
      \State $t_i = \left(\sum_{k = 1}^d \left(\max_{j \in B_i} x^{(p)}_{jk} - \min_{j \in B_i} x^{(p)}_{jk} \right)\right)^{1/2}$  \Comment{Sum of box dimensions surrounding $B_i$}
      \State $m = \argmin_{k \in B_i} \|x^{(p)}_k - X_i\|_*$  \Comment{Closest to target,$X_i$}
      \If{ $\min_{j \in B_i \cap \overline{C}^{(p)}} \|x^{(p)}_j - x^{(p)}_m\|_*
      \ge \tau_2\; t_i$} \Comment{No nearby $\overline{C}^{(p)}$ vertex}
        \State $\overline{C}^{(p)} = \overline{C}^{(p)} \cup \{m\}$
      \EndIf
    \EndFor
    \State $C^{(v)} = \mbox{\func{find\_colocated\_velocities}}(\overline{C}^{(p)})$
  \end{algorithmic}
\end{algorithm}

\RT{Algorithm \ref{a:midpoints}'s cost is clearly proportional to the number of
fine mesh vertices, $|F^{(p)}|$. Each fine vertex  (i.e, each $i$ in the {\bf
for} loop) requires a constant amount of work, though this constant is not
small.  Specifically, this consists of the computation of $X_i$ and the set
$B_i$ followed by calculations that are each proportional to the size of
$|B_i|$. The $X_i$ computation is proportional to the number of $C^{(p)}$
vertices that are used to interpolate to the $i^{th}$ fine vertex (at most $2^d$ in the ideal uniform case). The sizes of the $B_i$ are
bounded from above due to the PDE nature of the problem (meaning that for every
fine point there is a bounded number of coarse points that are close) and the
fact that coarse points interpolate to fine points only within distance three.
The $B_i$ determination requires visiting a single coarse point from $S_i$ (any
coarse point from the set is suitable) and checking if any of the fine points
that interpolate from it  also interpolate from all other points in $S_i$. On
regular grids with perfect distance-four coarsening (e.g.,
Figure~\ref{f:prescpts}), the number of visited fine points is $7^d-1$ and $|B_i|$
is bounded by $3^{d-1}7^{d-1}$ (where again $d$ is the problem dimension). }

The velocity prolongator sparsity pattern $N^{(v)}$ is defined exactly in
the same fashion as for $N^{(p)}$ \RT{with the exception that each velocity
component is treated separately. That is, a fine grid velocity associated
with a coordinate direction  only interpolates from coarse grid velocities
associated with the same coordinate direction (i.e., $P_\ell^{(v)}$ is 
block diagonal with each block corresponding to a different coordinate
direction).}
Specifically, nonzeros in the $i$th
row of $N^{(v)}$ are given by all $C^{(v)}$ vertices that are within
distance three of vertex $i$ within the nodal graph of $\widetilde{A}^{(v)}$.
% \fix{AP}{This is inconsistent with current interpretation of
% Algorithm~\ref{a:build_block_P}.}
This  $\widetilde{A}^{(v)}$ is defined by dropping
small entries from $A_{NS}$ as
previously discussed for $\widetilde{A}^{(p)}$.
%is defined by dropping
%indicated in the Appendix~\ref{s:appendix_dropping}.

\begin{figure}[t]
\centering
  \includegraphics[scale=.30]{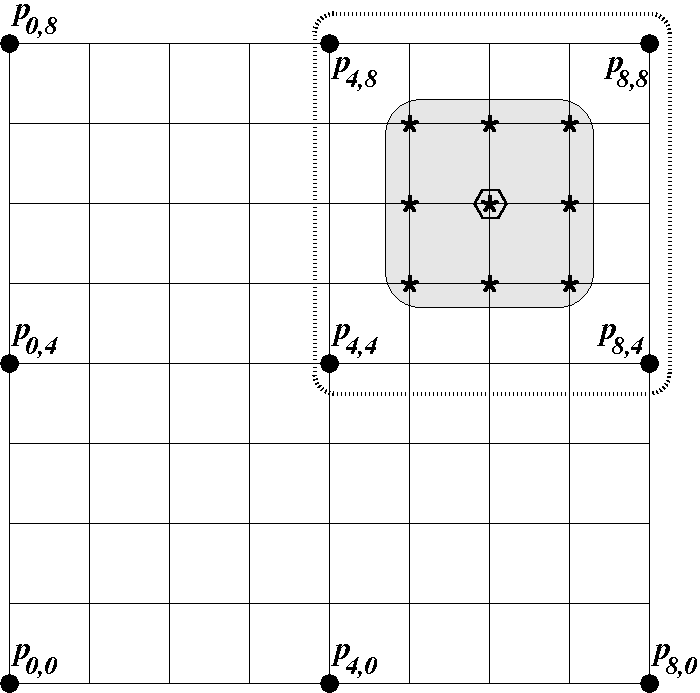}~~~~~~~
  \includegraphics[scale=.30]{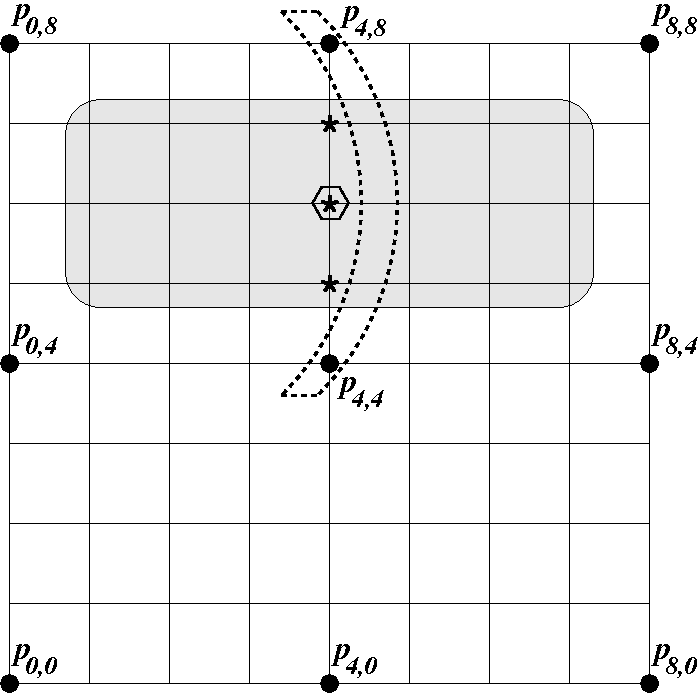}
\caption{Two mid-point scenarios.
Each `*' vertex interpolates from the same subset of $C^{(p)}$,
the solid dots within the dashed region.
The hexagon gives the target spatial location and
gray regions depict vertices that are searched to find
a vertex near the target. }
\label{f:midpoints}
\end{figure}

\begin{figure}[t]
\centering
  \includegraphics[scale=.30]{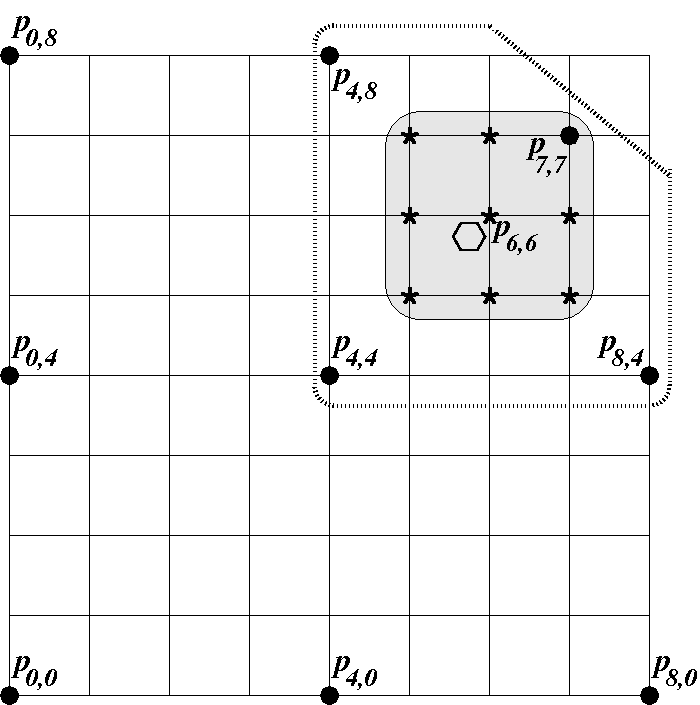}~~~~~~~
  \includegraphics[scale=.30]{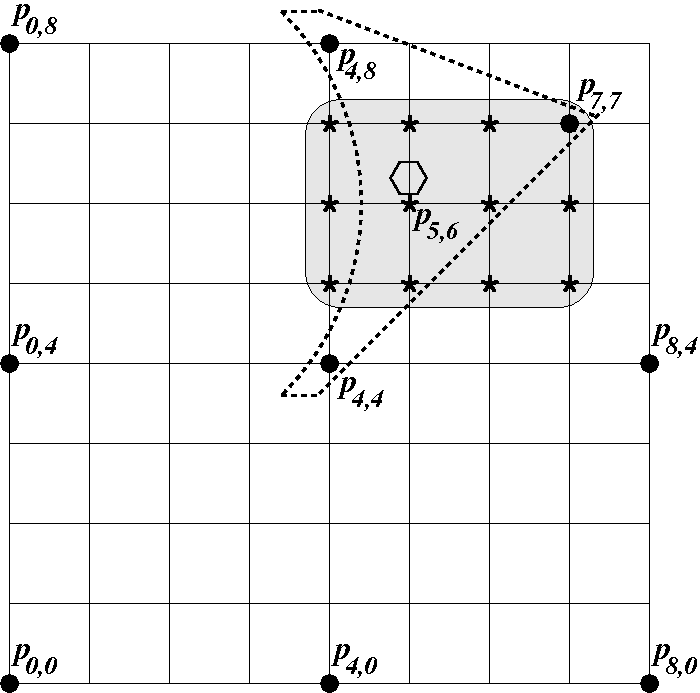}
\caption{Perturbed mid-point scenarios.
Each `*' vertex interpolates from the same subset of $C^{(p)}$,
the solid dots within the dashed region.
The hexagon gives the target spatial location and
gray regions depict vertices that are searched to find
a vertex near the target. }
\label{f:pertmidpoints}
\end{figure}

Figure~\ref{f:midpoints} illustrates an ideal case where
the original discretization is on a uniform regular mesh and
all $C^{(p)}$ vertices are equi-spaced. Specifically, the
right image\footnote{Again, the image
displays the underlying mesh as opposed to the matrix graph.}
depicts
 three vertices by a `*'. If we refer to these
vertices as $k, \ell, $ and $m$, then $S_k \equiv S_\ell \equiv S_m $
consists of two vertices in the banana shaped region.
%~\unsure{AP}{Do we really want "banana-shaped" region".}
The small hexagon gives
the target spatial location ($X_k \equiv X_\ell \equiv X_m $)
and vertices surrounded by the gray region define the
$B_k \equiv B_\ell \equiv B_m $, all within distance three from both
$S_k / S_\ell /S_m $ vertices.  In this situation the target location
is exactly at the same location as a pressure vertex. The left
image shows the same information for a different mesh location.
%where the associated $S$ sets include four vertices.
Figure~\ref{f:pertmidpoints} illustrates perturbed scenarios where
the $p_{7,7}$  vertex is in $\overline{C}^{(p)}$ as opposed to
the $p_{8,8}$  vertex. Thus, a heuristic has chosen a coarse point
that is a distance three (instead of four) from existing coarse
points. In this case, $p_{7,7}$ also lies within the $S$ sets associated
with the `*' vertices and the barycenters are offset.
If $p_{6,6}$ is first added to $\overline{C}^{(p)}$
(due to the left image computation), then after $p_{5,6}$ might
be added to $\overline{C}^{(p)}$ (the right image calculation)
or deemed too close to $p_{6,6}$
depending on the value of a user provided $\tau_2$.
That is, the total number of mid-points might vary
somewhat with either several nearby mid-points or a sparser
collection of mid-points. It is important to notice that
even if the mid-point selection has some irregularities,
this does not directly propagate to coarser levels. That is,
coarse velocities at any given level are essentially a function
of coarse pressures on that level and not a function
of coarse velocities from a preceding level. Thus, a less
than ideal mid-point selection does not have a direct
influence on the coarsening of velocities on even
coarser levels.
% \info{AP}{Should be careful here (see discussion associated with
% Algorithm~\ref{a:build_block_P}. So, in real life coarse velocities {\bf do}
% depend on the choice of fine level velocities. }
Finally, \RT{additional heuristics are used to convert some
$F^{(v)}$ vertices to $C^{(v)}$ vertices if
some $F^{(v)}$ vertices are distant from all current  $C^{(v)}$ vertices}.

\subsection{Determination of transfer coefficients using EMIN-AMG}

There are several ways in which interpolation coefficients can be
determined.  Here, we use the energy-minimizing framework,
EMIN-AMG, proposed in~\cite{Olson2011}, which is a generalization of ideas
in~\cite{Mandel1999} and closely related to \cite{Brandt00generalhighly,Brandt01multiscalescientific, BrZi2007, KoVa2006, Va2010, Wa2000, WaChSm2000,XuZi2004}.
We summarize only key features % of the energy-minimization AMG algorithm
as the details are not crucial to this paper. Flexibility is the most
significant aspect of EMIN-AMG for the present context. In particular,
the EMIN-AMG algorithm does not make any implicit or explicit assumptions
about the selection of coarse points or grid transfer sparsity patterns.
This frees one to choose these components in a way to preserve
important features. These features might be discretization
features (as in this paper) or they might be application features,
such as cracks or interfaces that one wishes to maintain throughout the
multigrid hierarchy. This flexibility is in contrast to popular AMG
methods such as classical AMG~\cite{BrMcRu1984,Ruge1987} or smoothed aggregation
AMG~\cite{Vanek1996} where the choice of coarse points (or aggregates),
interpolation sparsity pattern, and interpolation coefficients
is somewhat intermingled. Another AMG algorithm with similar
flexibility is given in~\cite{bootstrap}.

% The main idea with energy minimization AMG is ...

The main idea of EMIN-AMG is as follows. Let ${\cal N}$ be a set of
matrices with a specific (previously specified) nonzero pattern and dimensions,
and $W$ be a set of fine level modes (e.g., vectors) requiring exact
interpolation. Prolongator coefficients are determined through an
approximate solution of a constrained minimization problem
\begin{equation}
  P = \argmin\limits_P \sum\limits_j \|P_j\|_{\chi}^2
  \label{eq:min_problem} \\
%\end{equation}
~~~~~~\mbox{subject to}~ \\
%\begin{equation}
  P \in {\cal N},\quad\mbox{and}\quad W \in range(P) .
\end{equation}
Here, $\chi$ is some matrix norm,
$P_j$ is the $j$th column of $P$, and the sum is over all columns in $P$.
For symmetric problems (e.g., $A_S$ and approximations to the pressure
Schur complement), it is natural to take $\chi$ to be the $A$-norm.
Likewise, it is reasonable to take $\chi$ to be the $A^T A$-norm
for non-symmetric systems.
While the solution of an optimization problem
may seem onerous compared to the task of solving a linear system,
only a rough approximate solution to \eqref{eq:min_problem} is needed
and an iterative process can be started with a simple but good initial
prolongator satisfying the constraints.

In our case, $W$ is just
a single vector of all ones when constructing $P^{(p)}$. This corresponds
to the requirement that a constant can be perfectly recovered by interpolation.
For $P^{(v)}$, $W$ consists of two (three) vectors in two (three) dimensions.
Each vector contains only zeros and ones and corresponds to a constant
for a velocity component in each of the different coordinate directions.
It follows that the constraint $ W \in range(P) $ can be satisfied if the
sum of all nonzeros in each row of $P$ is one (assuming that
the sparsity pattern
does not mix different velocity direction components). That is,
a coarse level constant interpolates to a fine level constant. Thus,
a simple initial feasible guess for an optimization algorithm applied
to \eqref{eq:min_problem} just takes the binary sparsity pattern matrix and
divides each entry by the number of nonzeros in each row. % of $N$.
The minimization problem is solved iteratively, for instance, using
Algorithms 2 or 3 from~\cite{Olson2011}. These algorithms essentially
correspond to applying a
constrained version of CG or GMRES to solve for a
prolongator matrix. The key is that the overall cost of the prolongator
construction is only modestly higher than that of a more conventional
AMG procedure as only a couple of EMIN-AMG iterations are sufficient.

In symmetric cases, $R$ is taken as $P^T$. However,
for the Navier-Stokes problem a Petrov-Galerkin projection is employed, i.e.
$R \ne P^T$. An energy minimization procedure is also used for determining
the $R$ matrix. This procedure employs the transposed sparsity pattern for $P$
and uses the $A A^T$-norm for $\chi$. This choice of $R$ generally
leads to projected discretization operators with satisfactory stability
properties even in the presence of strong convection. This is discussed
further in~\cite{Olson2011,WiTuWaGe2012}

\section{Smoothers}\label{s:smoothers}
Our choice of smoothers is somewhat restricted as typical AMG smoothers
such as Jacobi and Gauss-Seidel are ineffective on saddle-point systems
due to their negative eigenvalues.  Relatively standard incomplete
factorizations such as ILU(1) can be used for AMG smoothers. As the sparsity
pattern of the incomplete factors is closely connected to the sparsity pattern
of the original matrix (e.g. the ILU(0) factors have the same sparsity pattern
as the initial discretization matrix), the zero block in the discretization
matrix can cause issues.  In this paper, our ILU(1) implementation treats the
matrix diagonal (including those within the zero block) as being nonzero to
encourage fill-in within the part of the incomplete factors associated with
the zero block. In addition to ILU(1), we consider two families of smoothers
that specifically take advantage of the block structure: Vanka and
Braess-Sarazin relaxation.

\subsection{Vanka relaxation}
Vanka smoothing was originally proposed in~\cite{Vanka1986} for
finite-difference schemes. Further analysis of Vanka methods for finite-element
discretizations of the Stokes equations has been done in~\cite{MacLachlan2011}.
The Vanka scheme corresponds to an overlapping block Gauss-Seidel method.
The blocks are defined by partitioning all dofs into overlapping sets
$T_i$, $i = 1, \dots, n^{(p)}$.
The number of sets is the same as the number of pressure dofs when \qq elements
are employed.
Each $T_i$ can be defined algebraically by taking all column indices
corresponding to nonzero entries in the $i$th row of $B$ along with the
$n^{(v)}+i$th index.  That is, each $T_i$ consists of a single pressure
dof and all velocity dofs that are either co-located or adjacent (in the
matrix graph) to this single pressure dof.
This choice of blocks is motivated by the saddle point nature of the
problem so that
each block submatrix is also a saddle point matrix.
To apply one step of Vanka relaxation to a linear
system, $A{\bf x} = {\bf f} $ where the current approximation
is given by ${\bf x^k} $, one computes an update of the form
\begin{equation}
  {\bf x^{k+1}} = {\bf x^k} + \omega R_i^T (R_i A R_i^T)^{-1} R_i \left({\bf
  f} - \mathcal{A}{\bf x^k}\right).
\end{equation}
Here, $R_i$ is a binary projection operator restricting a global vector to a local
one corresponding to dofs in $T_i$, and $\omega$ is the under-relaxation
parameter.  The overall relaxation procedure is done in
a Gauss-Seidel manner cycling through all sets $T_i$.

\subsection{Braess-Sarazin relaxation}
Braess-Sarazin algorithms were originally considered as a relaxation scheme
for Stokes problems~\cite{Braess1997,Larin2008}, and later applied to
Navier-Stokes equations~\cite{John2000}. Compared to Vanka relaxation, where
the smoothing procedure relies on solving multiple local saddle-point problems,
Braess-Sarazin relaxation solves a global problem, though greatly simplified.

A single step of the relaxation procedure can be written as
\begin{equation}
\begin{pmatrix} \bf{u^{k+1}} \\ \bf{p}^{k+1} \end{pmatrix} =
\begin{pmatrix} \bf{u^k} \\ \bf{p}^k \end{pmatrix} + \begin{pmatrix}\frac{1}{\omega}
D & B^T \\ B & O \end{pmatrix}^{-1}\left(
\begin{pmatrix} \bf{f} \\ \bf{g} \end{pmatrix} -
  \begin{pmatrix}
    \widehat{A} & B^T \\ B & O
  \end{pmatrix}\begin{pmatrix} \bf{u^k} \\ \bf{p}^k \end{pmatrix}
  \right),
\label{eq:bs}
\end{equation}
where $\widehat{A}$ is the velocity block of $A$, $D$ is a suitable
preconditioner for $\widehat{A}$, and $w$ is a relaxation parameter. As
\begin{equation}
\begin{pmatrix}\frac{1}{\omega} D & B^T \\ B & O \end{pmatrix}^{-1} =
\begin{pmatrix}I & \frac{1}{\omega} D^{-1}B^T \\ O & I \end{pmatrix}^{-1} \begin{pmatrix}
\omega D & O \\ B & -\frac{1}{\omega} BD^{-1}B^T \end{pmatrix}^{-1},
\label{eq:bs2}
\end{equation}
a solution of a Schur complement system, $S = BD^{-1}B^T$, is required. In our experiments, the Schur complement is explicitly formed and solved
approximately by a relaxation procedure, e.g. via %by performing
several Gauss-Seidel iterations.
The Braess-Sarazin smoother requires a practical choice of the matrix $D$.
The original paper considered $D = \diag(\widehat{A})$. While that choice
performed reasonably well in our experiments, we found that faster convergence
is achieved \RT{with $D$ being a diagonal matrix with entries $d_{ii} = \sum_j
|(\widehat{A})_{ij}|$. That is, $d_{ii}$ is the sum of the absolute values in the $i$th row
of $\widehat{A}$}.
%\improve{AP}{Why? Is that true? Is there some data supporting this? Do
%we need to investigate?}
%\improve{AP}{SIMPLE-C related?}

\section{Numerical results}\label{s:results}
In this section, Stokes and Navier-Stokes problems are considered.
$\bm{Q}_2$ elements are used to discretize velocities and $\bm{Q}_1$
elements are used to discretize pressure.
The discrete problems were generated using the
{\sc \textsf{IFISS}} software package~\cite{IFISS} (version 3.3) \RT{written in
MATLAB. The proposed algorithms were prototyped in MATLAB using the {\sc
\textsf{MueMat}} package~\cite{MueMat}, and later implemented in C++ in
the {\sc \textsf{MueLu}} multigrid package~\cite{MueLu}. All numerical results were
produced with {\sc \textsf{MueLu}} (as of Trilinos version 12.6) with
a single exception of the unstructured circle-driven cavity problem that was
solved with {\sc \textsf{MueMat}}}.

In all of our calculations,  five Gauss-Seidel iterations are performed
to approximate the solution of the Schur complement system within the
Braess-Sarazin smoother (BS), and the
relaxation parameter is set to $0.666$. The Vanka smoother's under-relaxation
parameter is fixed at $0.5$, though it was observed in~\cite{Vanka1986} that
higher Reynolds numbers would benefit from lower parameter values. Unless stated
otherwise, GMRES is
used as an iterative method with the residual
tolerance for the stopping criteria chosen to be $10^{-6}$.
All results in this section use $\tau_1 = 0.06$ in~\eqref{e:filtering} for dropping small entries while
constructing $\widetilde{A}^{(v)}$ and $\widetilde{A}^{(p)}$, and
$\tau_2 = \sqrt{1.5\cdot 10^{-3}}$ in Algorithm~\ref{a:midpoints} to decide
whether a candidate coarse velocity mid-point is sufficiently far from already
chosen velocity coarse points so that it should be added as an additional coarse
point.

\subsection{Stokes lid-driven cavity problem}

\begin{figure}[t]
\centering
  \includegraphics[scale=0.30]{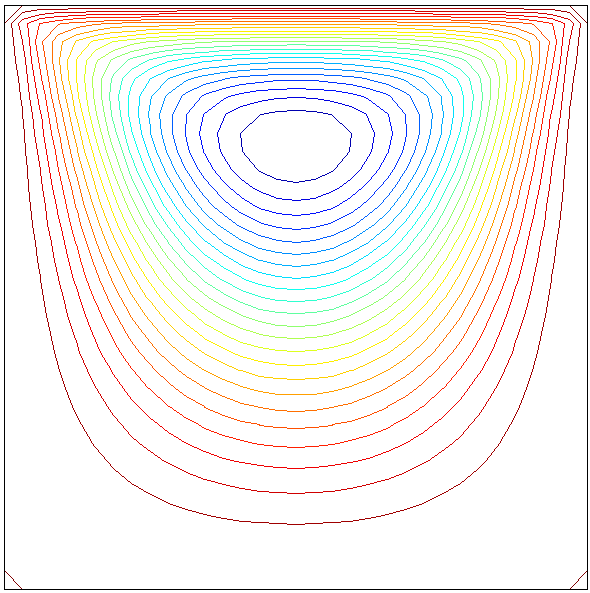}
\caption{Equally distributed streamlines for the Stokes lid-cavity problem.}
\label{f:stokes_lid_cavity}
\end{figure}

We begin with a benchmark lid-driven cavity Stokes problem
%classic fluid dynamics problem, known as
%The Stokes problem is considered
on a square $\Omega = (-1,1)^2$ domain
with leaky boundary conditions. Specifically, the $y$ component of velocity
is zero on the boundary, i.e. $u^{(y)}|_{\partial\Omega} = 0$. The
$x$ component of velocity is only nonzero on the top of the cavity.
In particular,
\begin{equation}
  u^{(x)} = \begin{cases}
    1,& \mbox{for}\quad y=1; -1 \le x \le 1, \\
    0,& \mbox{otherwise}.
  \end{cases}
\end{equation}
We consider a set of uniform meshes. The discrete
system has one zero eigenvalue corresponding to constant pressure.
A typical solution on a $10 \times 10$ mesh is illustrated in
Figure~\ref{f:stokes_lid_cavity}. Prolongator coefficients are computed
with a single EMIN-AMG CG step~\cite{Olson2011}. One iteration of the Vanka
smoother or two iterations of Braess-Sarazin smoother are considered as pre-
and post-smoothers in the multigrid hierarchy. Vanka relaxation is used on
the coarsest grid instead of a direct solver to avoid direct solver issues
associated with the singularity of the coarsest grid problem.

\begin{table}[t]
\caption{Number of iterations and multigrid run times for the Stokes lid-driven cavity.}
\label{t:stokes_lid_cavity}
\centering
\tabsize
\begin{tabular}{ccccrrcrrcrr}
  \toprule
  \multirow{2}{*}{Dofs} & \multirow{2}{*}{Complexity} &
  \multirow{2}{*}{$\ell_{max}$} & \multicolumn{3}{c}{Vanka (1,1)} &
  \multicolumn{3}{c}{BS (2,2)} & \multicolumn{3}{c}{ILU(1) (1,1)}\\
  \cmidrule(lr){4-6} \cmidrule(lr){7-9} \cmidrule(lr){10-12}
  &&& Its & Setup & Solve & Its & Setup & Solve & Its & Setup & Solve \\
  \midrule
     659 & 1.08 & 2  & 11 &  0.02 &    0.06 &  9 &  0.01 &   0.02  &  6 &  0.02 & 0.01 \\
    2467 & 1.08 & 2  & 12 &  0.05 &    0.20 & 10 &  0.03 &   0.03  &  7 &  0.05 & 0.02 \\
    9539 & 1.08 & 3  & 14 &  0.21 &    1.02 & 14 &  0.09 &   0.10  &  9 &  0.19 & 0.06 \\
   37507 & 1.08 & 3  & 15 &  0.93 &    4.38 & 19 &  0.36 &   0.36  & 10 &  0.77 & 0.22 \\
  148739 & 1.11 & 4  & 19 &  5.11 &   22.31 & 27 &  2.17 &   2.19  & 13 &  4.08 & 1.16 \\
  592387 & 1.05 & 5  & 18 & 28.84 &   84.80 & 20 & 13.34 &   6.70  & 17 & 20.32 & 5.96 \\
  \bottomrule
\end{tabular}
\end{table}

Table~\ref{t:stokes_lid_cavity} summarizes the results, including number of
iterations and run times using a multigrid preconditioned \RT{GMRES} solver. Our
aggressive approach to coarsening (choosing
coarse pressure dofs that are distance four from each other) results
in a modest number of multigrid levels, where the
coarsest hierarchy level for all examples has less than 205 total dofs.
Further, the multigrid operator complexities are very small.
These operator complexities are defined as the sum of the
number of nonzeros of all discretization matrices on all levels divided
by the number of nonzeros for the finest level discretization matrix.
As the table illustrates, the number of additional nonzeros associated with the coarse
level matrices is quite small and one can expect that the storage and
computational time associated with these coarse operators is also quite modest.
That is, the cost per iteration does not grow appreciably when more levels are
employed.
It is also particularly important that coarse discretization stencil widths do not
grow too large on large scale parallel machines as large stencil widths imply
longer distance communication.

Overall, the number of iterations remains relatively stable with Vanka and Braess-Sarazin
smoothing which suggests an $h$-independence property of our approach (i.e., the number
of iterations remains bounded as the mesh is refined). The generally well-behaved
nature of the convergence rates is an indication that the AMG hierarchy (generated
by the proposed coarsening, sparsity patterns, and EMIN-AMG) is suitable for this
Stokes problem.  Though the iterations
do grow slowly with ILU(1) relaxation, it does require the least time in the solve
phase for the meshes considered in this study. It is well known that the ordering
of unknowns can have a significant influence on the convergence behavior when ILU methods
are employed~\cite{Vuik2008,Dahl1992}. For these Stokes experiments we found that a natural
ordering (i.e., lexicographical starting from the lower-left corner) worked
fine.  However, a reverse Cuthill-McKee (RCM) ordering is needed to obtain satisfactory
convergence rates with ILU for all the remaining Navier-Stokes examples shown
in this paper.
%Overall, AMG with Vanka smoothing requires slightly less  faster convergence
%rate. This may be due to its nature of preserving global properties in the
%local sub-problems.
Overall, the Vanka run times are noticeably slower than
Braess-Sarazin run times even though the convergence rate is a bit better.
This iteration/run time behavior between Vanka
and Braess-Sarazin has also been observed in~\cite{TomPaper}. The Vanka smoother
involves numerous dense matrix solves
associated with each Vanka block, which are of larger size due to the second
order approximation in velocities. Additionally, residual components
associated with velocity are updated multiple times by the Vanka smoother
due to a significant overlap among Vanka blocks. It should be pointed out
that our smoothing code has not been fully optimized.

With our current implementation, the setup phase requires more time than the
solve phase.  In general, the initialization time associated with the Braess-Sarazin
smoother is small relative to the total AMG setup time while the Vanka and ILU(1) setup
times are large. Roughly, the Vanka smoother setup time (or ILU(1) setup time) can
be estimated by taking the difference between the total AMG setup time using Vanka (or ILU)
and the total AMG setup time using Braess-Sarazin. Even with the Braess-Sarazin smoothing,
the setup time is noticeable. The overwhelming majority of this time is spent on graph
distance calculations, used in the coarsening and sparsity pattern calculation, and not
in the EMIN-AMG algorithm. We believe that the large time for distance calculations is due to
a poor implementation as opposed to intrinsic to the algorithm. While we have not
further optimized this graph phase of the algorithm, it is important to note that the graph
calculation is easily amortized over a nonlinear sequence of linear systems. That is, the graph
coarsening can be performed just once over a sequence associated with a nonlinear solve. For the
finest mesh associated with Table~\ref{t:stokes_lid_cavity}, the graph setup
phase required 9.74 seconds. As this cost only needs to be incurred once, the
incremental setup cost for an additional linear system using Braess-Sarazin
would only be 3.60 seconds, which is about half the time of the solve phase.
%\improve{RST}{Not sure if we should say more or less? We need the XYZ times}
%\improve{AP}{I'm OK with the current phrasing.}
%This is partially due to the
%The setup of the Vanka smoother involves a large number of
%dense matrix factorizations
%On the other
%hand, the construction of the Braess-Sarazin smoother involves a single matrix
%product $BD^{-1}B^T$. Similarly, the
%while Braess-Sarazin updates those components only
%once performing iterations only on the pressure components.

\RT{
Let us now examine the effect of the parameters $\tau_1$ and $\tau_2$ on our
results. The parameter $\tau_1$ was introduced in Section~\ref{s:ptransfer}. It
is used for stencil reduction in matrices associated with velocity and pressure
blocks. From a convergence point of view, smaller values of $\tau_1$ correspond
to denser block matrices used in prolongators construction and thus result in
more aggressive coarsening. On the other hand, larger values of $\tau_1$
correspond to sparser matrices and poorer approximations to the original
matrices. Table~\ref{t:tau1} demonstrates the
convergence for the lidcavity problem on a fixed mesh with varying values of
$\tau_1$ while $\tau_2$ is fixed. The number of iterations reaches minimum at
$\tau_1 = 0.05$ and increases thereafter. The size of the first coarse matrix
increases from 28088 rows for $\tau_1 = 0$ to 74428 for $\tau_1 = 0.25$.
Therefore, the selected $\tau_1$ is indeed a suitable choice for the problem.
We expect it to be reasonable for Navier-Stokes problems with low Reynolds numbers.
}

\RT{
The parameter~$\tau_2$ was described in Section~\ref{s:vtransfer}. It is responsible for
determining the distance at which close points are considered to coincide. We
found that our results are insensitive to this parameter. For instance, we
found that the number of iterations remained constant for varying $\tau_2$
from $0.0$ to $0.4$. The parameter produced only a minor difference in the size of
some coarse level matrices (e.g., it slightly reduced the size of the
level 2 matrix for the 592387 mesh from 2739 for $\tau_1 = 0.01$ to 2199 for $\tau_2 = 0.4$).
}

\begin{table}[t]
\caption{Number of iterations for the Stokes matrix of size 592387 with
  Braess-Sarazin smoother for varying values of $\tau_1$.}
\label{t:tau1}
\centering
\tabsize
\begin{tabular}{cc}
  \toprule
      $\tau_1$ & Number of iterations \\
    \midrule
      0.00     & 100+  \\
      0.05     & 20    \\
      0.10     & 21    \\
      0.15     & 32    \\
      0.20     & 36    \\
      0.25     & 100+  \\
  \bottomrule
\end{tabular}
\end{table}

\RT{
We investigate the stability of the coarse grid operators by using a technique
suggested in~\cite{boffi13} (section 3.4.3) where a relationship between the
smallest non-zero singular value of an off-diagonal matrix (in our case the scaled divergence operator) and the inf-sup constant
$\beta$ is established. Instead of a looking at a sequence of matrices based on mesh
refinement, we consider a sequence matrices derived from coarse level operators.
Specifically, let $M_0^{(p)}$ and $M_0^{(v)}$ be pressure and velocity mass
matrices associated with the fine level discretization. We denote by
$M_{\ell}^{(p)}$ and $M_{\ell}^{(v)}$ the projections of these matrices to
coarse levels, i.e.
\begin{equation}
  M_{\ell}^{(p)} = R_{\ell}^{(p)} M_{\ell-1}^{(p)} P_{\ell}^{(p)}, \quad
  M_{\ell}^{(v)} = R_{\ell}^{(v)} M_{\ell-1}^{(v)} P_{\ell}^{(v)}, \quad
  \ell = 1, \dots, \ell_{max}.
\end{equation}
Let us also denote by $B_{\ell}$ the (2,1) block of the coarse operator
$A_{\ell}$, corresponding to pressure rows and column velocities. Clearly,
$B_0 = B$ from equation~\eqref{eq:stokes_block_system}. We then
compute the smallest non-zero singular values $\sigma_{min,\ell}$ of the
matrices
\begin{equation}
  \widetilde{B}_{\ell} =
  \lump(|M_{\ell}^{(p)}|)^{-\frac12} B_{\ell}
  \lump(|M_{\ell}^{(v)}|)^{-\frac12}.
\end{equation}
}
\begin{table}[t]
\caption{Smallest non-zero singular values of matrices $\widetilde{B}_{\ell}$ for the Stokes
  matrix of size 592387.}
\label{t:infsup_estimate}
\centering
\tabsize
\begin{tabular}{ccccc}
  \toprule
    $\ell$              & 0    & 1    & 2    & 3    \\
    \midrule
    $\sigma_{min,\ell}$ & 1.54 & 1.56 & 1.56 & 1.07 \\
  \bottomrule
\end{tabular}
\end{table}
\RT{
The results are shown in Table~\ref{t:infsup_estimate}.
%\footnote{The coarsest
%level is not shown due to having a single pressure degree of freedom resulting
%in $B_5 = O$.}.
The smallest singular values are clearly separated from 0 which suggests a
similar stability of the operators, though
%While the values for
the coarsest operator's minimum singular value
is somewhat lower (due in part to the small size, $150$ pressure
dofs, of this matrix).
%of size less than 150).
%and could be considered to correspond to a very rough mesh.
The singular values
for finer level matrices are all very close to each other.
% In this example, this could be explained by the fact that our algorithm
% automatically constructs coarsening which is similar to a standard geometric
% coarsening.
}

\subsection{Stokes circle-driven cavity problem}
\RT{
%% Adaptation of the ECE 100 template located at
%% https://www.overleaf.com/latex/templates/ece-100-template/pjrrfybfggqt#.Vg2ScOlqnky
%% Thanks to author Patrick Bartman, permission via Creative Commons CC by 4.0
%%
%\documentclass[a4paper, 11pt]{article}
%\usepackage{lipsum} %This package just generates Lorem Ipsum filler text.
%\usepackage{fullpage}
%\usepackage{graphicx}
%\usepackage{color}
%\usepackage{booktabs}
%\usepackage{multirow}
%%\usepackage{endfloat}
%\usepackage[compact]{titlesec}
%\titlespacing{\section}{0pt}{4pt}{0pt}
%\titlespacing{\subsection}{0pt}{2pt}{0pt}
%\titlespacing{\subsubsection}{0pt}{2pt}{0pt}
%
%\usepackage[normalem]{ulem}
%\newcommand{\irina}[2]{{\color{red}{#1}}
%{\color{red}{\sout{#2}}}}
%
%
%\begin{document}
%%Header-Make sure you update this information!!!!
%\noindent
%\begin{center}
%\textbf{Small addition on running mixed Q2-Q1 solver on unstructured meshes}
%\end{center}

To test the solution strategy on unstructured meshes, we consider
a modified form of the Stokes lid-driven cavity problem which we will call a
circle-driven cavity problem. Specifically,
the problem domain now corresponds to a square with a hole removed from
the center. Figure~\ref{sample mesh} illustrates one of the coarser meshes.
This mesh was generated with the Cubit software package~\cite{cubit}. While the
mesh resembles a structured mesh near the domain corners, it is unstructured
near the circular boundary. For our experiments, both velocity components
are set to zero (via Dirichlet boundary conditions) along the box boundary.
Additionally, Dirichlet conditions are applied on the circle boundary such
that the velocity along the circle has a magnitude of one and is oriented
in the clockwise tangent direction. That is, the horizontal (vertical)
velocity is $1$ at the topmost (leftmost) side of the circle and $-1$
at the bottom (rightmost) side of the circle.
Figure~\ref{sample solution} illustrates the computed solutions on the mesh associated with 7097 elements. The dark ring in the figure corresponds to the
circle boundary for $z = 0$ (and is drawn to help clarify the plot).
\begin{figure}[t]
\centering
  \includegraphics[width=3in,height=2.5in]{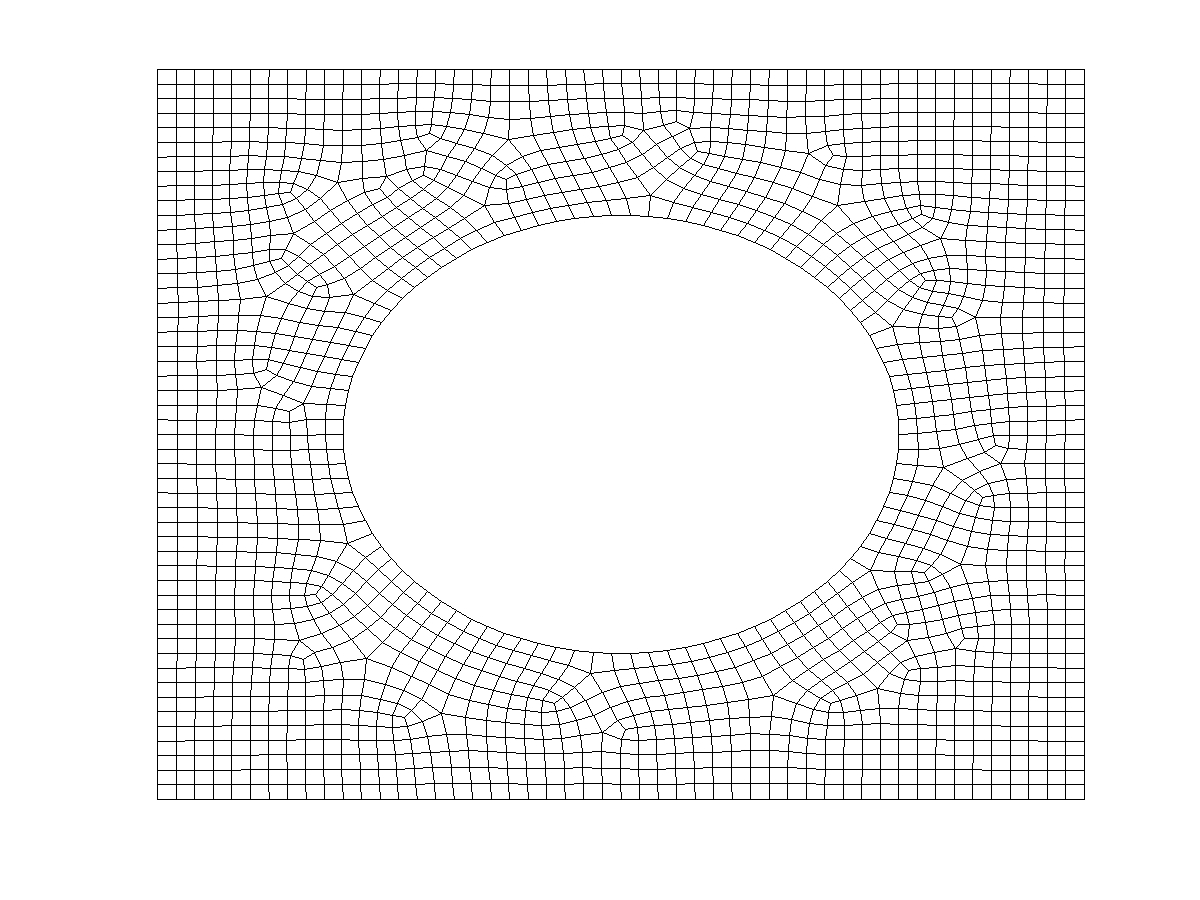}
\caption{Unstructured mesh for the Stokes circle-driven cavity problem containing 1756 quadrilateral elements.
\label{sample mesh}}
\end{figure}
\begin{figure}[t]
\centering
\includegraphics[trim=2.1in 2.0in 1.7in 2.0in, clip, scale=.20]{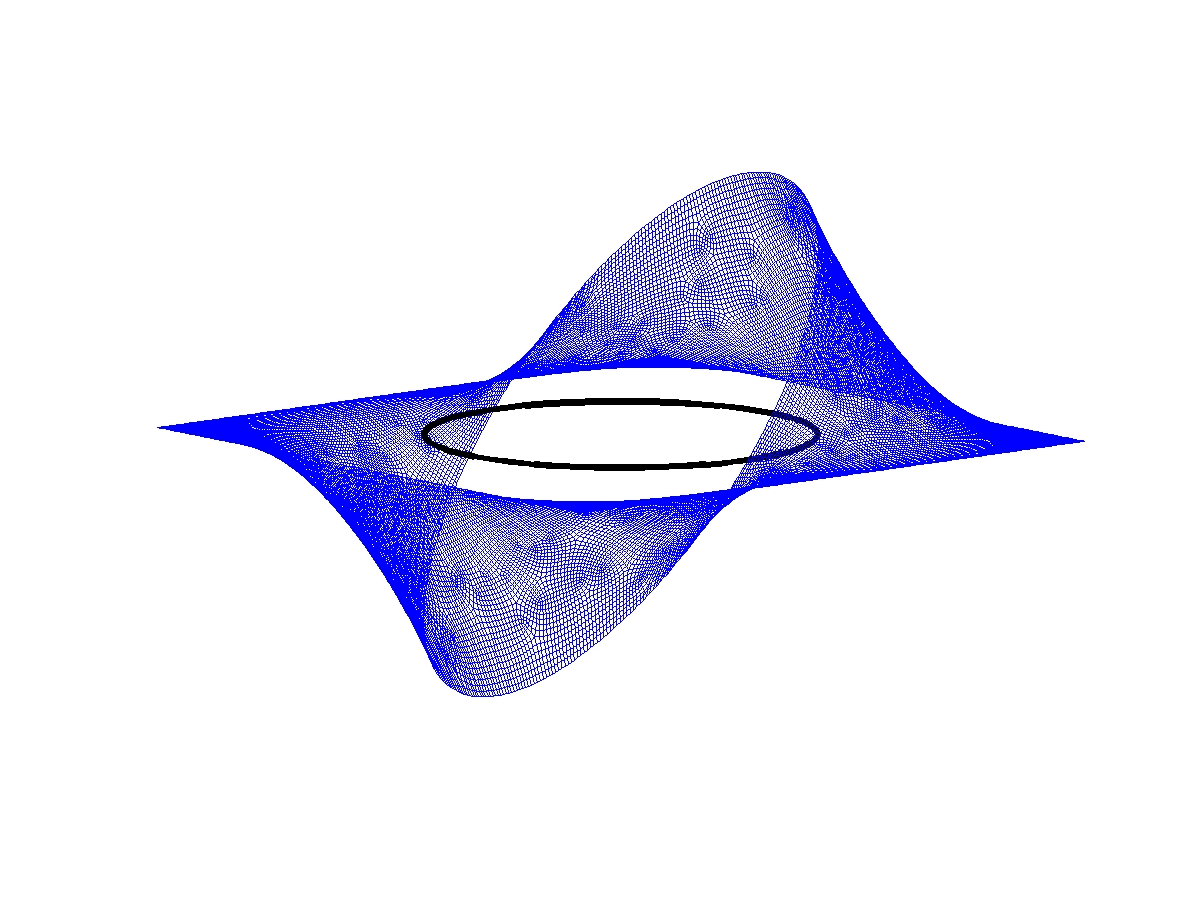}
\includegraphics[trim=2.1in 2.0in 1.7in 2.0in, clip, scale=.20]{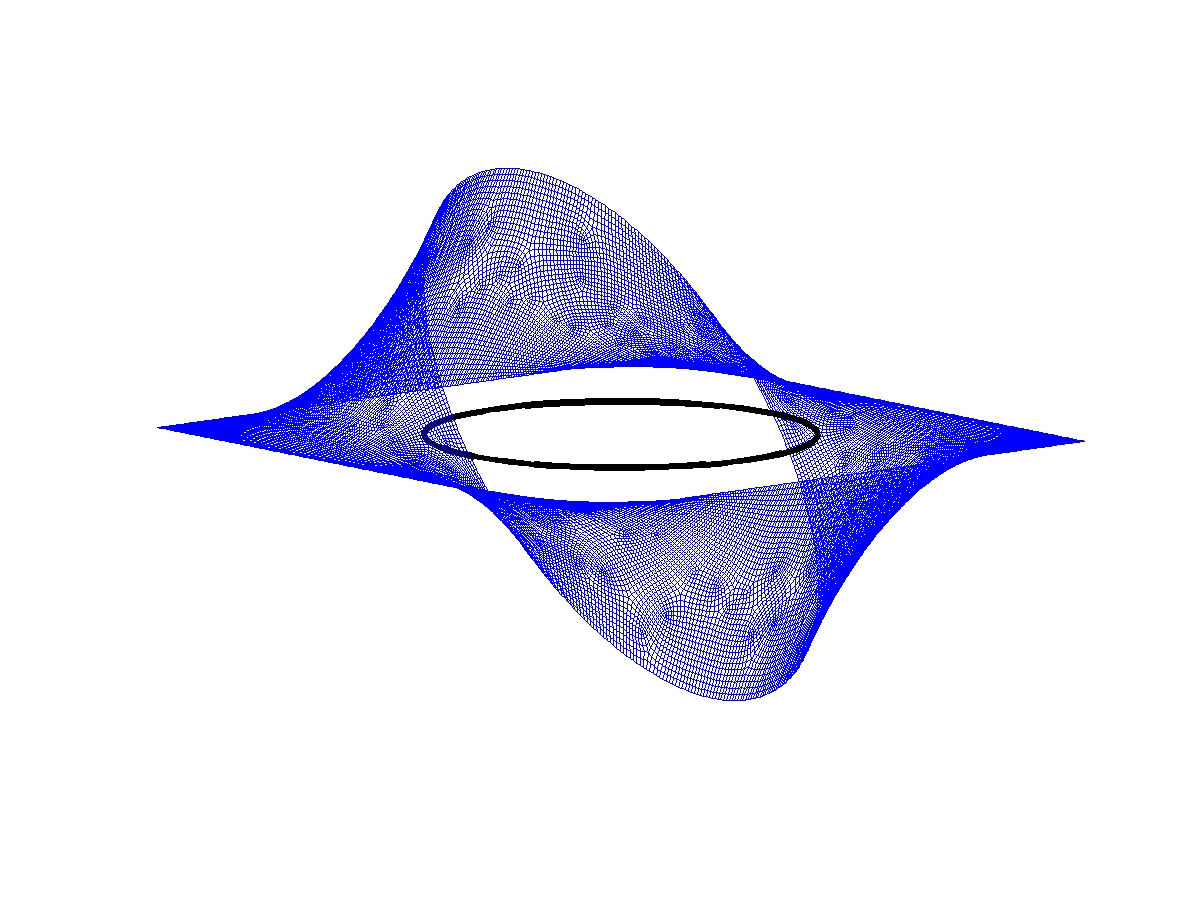}
  \caption{Computed horizontal (left) and vertical (right) velocity components
  for the Stokes circle-driven cavity problem. The rightmost portion of each figure
  corresponds to the top right corner of the cavity.
  \label{sample solution}}
\end{figure}
One can see that the horizontal and vertical prescribed values on the circular
boundary gradually decay as they approach the outer cavity boundary.

Table~\ref{t:unstructured} records the required number of GMRES iterations using the proposed multigrid preconditioner
with both Vanka and Braess-Sarazin smoothing.
\begin{table}[t]
\caption{Number of GMRES iterations for the Stokes circle-driven cavity problem.}
\label{t:unstructured}
\centering
\tabsize
\begin{tabular}{ccc|cc}
  \toprule
  Dofs & Complexity & $\ell_{max}$ & Vanka (1,1) & BS (2,2) \\
  \midrule
     1045 & 1.10 &  2 &  12 &   9  \\
     4078 & 1.12 &  3 &  13 &  12  \\
    16539 & 1.11 &  4 &  13 &  14  \\
    65343 & 1.10 &  4 &  14 &  18  \\
   261074 & 1.10 &  4 &  17 &  22  \\
  1043963 & 1.10 &  5 &  28 &  62  \\
  \bottomrule
\end{tabular}
\end{table}
One can see that the iteration count grows quite modestly with respect to
increases in mesh resolution (corresponding to a $1000$x increase in linear
system size) similar to the structured grid tests with a slight uptick for the
largest mesh. The exception is the Braess-Sarazin result on the largest mesh
which we believe to be a result of a deficiency in Braess-Sarazin smoothing on
steady problems. We have noticed the same sort of behavior in conjunction with
geometric multigrid for meshes with irregular mesh spacing~\cite{TomPaper}.

%\end{document}

}

\subsection{Navier-Stokes problems}

\begin{figure}[t]
\centering
  \includegraphics[scale=0.30]{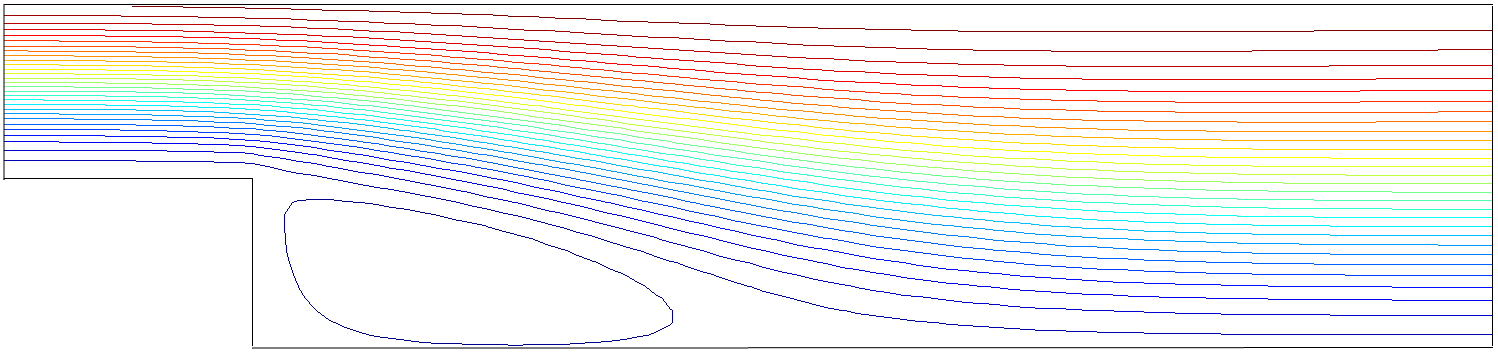}
\caption{Equally distributed streamlines for the Navier-Stokes backward-facing step problem.}
\label{f:ns_backward_step}
\end{figure}

\begin{figure}[t]
\centering
  \includegraphics[scale=0.30]{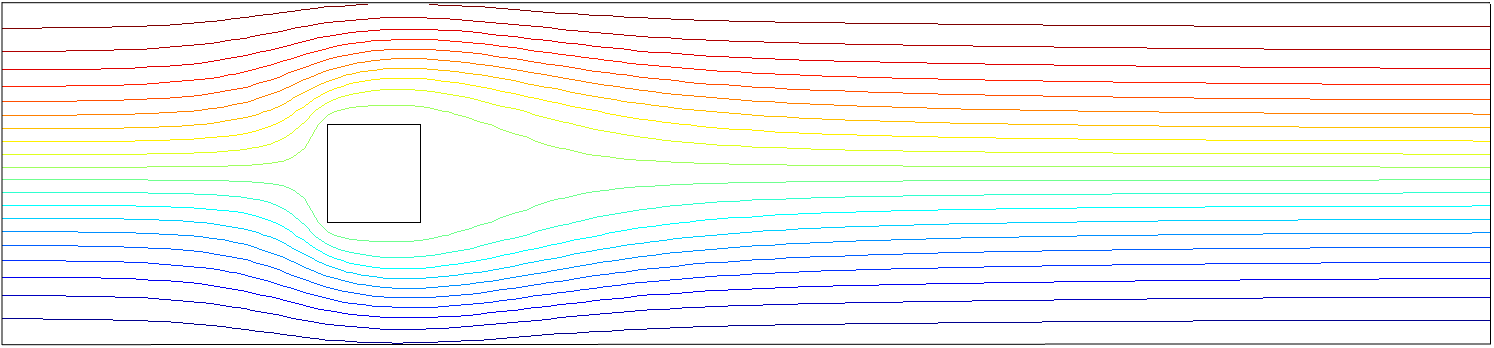}
\caption{Equally distributed streamlines for the Navier-Stokes obstacle problem.}
\label{f:ns_obstacle}
\end{figure}

Next, we consider the following problems for Navier-Stokes flow:
\begin{enumerate}
  \item Lid-driven cavity problem

    The lid-driven mesh and boundary conditions are identical to the Stokes
    setup. Two values of $\nu$ are considered, $\nu = 0.01$ and $\nu = 0.002$
    (corresponding to $\mathcal{R} = 200$ and 1000, respectively).

  \item Backward facing step problem

    An L-shaped domain is considered with a uniform mesh. A Poiseuille flow
    profile is imposed on the inflow boundary $(x=-1; 0 \le y \le 1)$, and a
    no-flow (zero velocity) condition is imposed on the walls. The outflow
    boundary $(x = L; -1 < y < 1)$ condition is set to Neumann, which
    automatically adjusts the mean outflow pressure to zero. Two values of $\nu$
    are considered, $\nu = 0.02$ and $\nu = 0.005$ (corresponding to
    $\mathcal{R} = 100$ and 400, respectively). The longer channel length $L = 10$
    is used for larger Reynolds numbers to allow for the exit flow to become
    well-developed (i.e., to have an essentially parabolic profile), while $L = 5$
    is used for smaller $\mathcal{R}$. Figure~\ref{f:ns_backward_step}
    demonstrates a typical solution for $\nu = 0.02$.

  \item Obstacle problem

    The obstacle problem has a similar setup to that of the backward facing
    step. In particular, $L = 8$ (for both $\nu = 0.02$ and $\nu =
    0.005$) and a uniform mesh is used to represent the domain. Additionally, a
    Poiseuille flow profile is imposed on the inflow boundary, no-flow condition
    is imposed on the walls and the same outflow boundary condition is imposed
    at the exit. Figure~\ref{f:ns_obstacle} demonstrates a typical solution for
    $\nu = 0.02$.
\end{enumerate}
As all of these problems are nonlinear, the matrix for AMG testing was chosen to be
the last linear system within a converged Picard sequence that is terminated when
a 2-norm of a nonlinear residual is less than $10^{-8}$.
Prolongator coefficients and restriction coefficients are each computed with a single
EMIN-AMG GMRES step.
One iteration of the Vanka smoother as well as two iterations of Braess-Sarazin
smoother are considered as a pre- and post-smoothers in the multigrid
hierarchy.

As the relative timing behavior is similar to the Stokes problem, we focus
on convergence results for the Navier-Stokes problems. These
are summarized in Table~\ref{t:ns_backward_step},
Table~\ref{t:ns_lidcavity}, and Table~\ref{t:ns_obstacle} for each of the
studied problems. Overall, the most important aspect of the convergence
behavior is that the number of iterations is relatively stable with
respect to mesh refinement, though there is certainly variation across
the methods/problems and occasionally some modest growth. This generally
gives us a certain level of confidence in the grid transfers and the overall
stability of the coarse grid operators produced by the proposed combination
of coarsening, grid transfer sparsity patterns and EMIN-AMG.

Similar to the Stokes problem, ILU smoothing generally leads to relatively
consistent growth in iterations as the mesh is refined, though iterations
dropped for the finest mesh lid-driven cavity problem when $\nu = 0.005$.
Some combinations of $\nu = 0.005$ and coarse meshes are problematic with ILU smoothing.
Though we have not performed a detailed study of these cases, we do note that
for these mesh sizes there are stability concerns for the finest level Galerkin
discretization when $\nu = 0.005$.  As with the Stokes problem, ILU smoothing
often leads to the fewest iterations, though Vanka smoothing leads to the most scalable
method with respect to iterations and mesh refinement. Interestingly, the number of
iterations drops steadily when AMG is used with Vanka smoothing for the $\nu = 0.005$
obstacle, though in other situations iterations are either nearly flat or rise modestly.
While convergence rates are relatively constant with respect to mesh refinement,
there is sensitivity to the Reynolds number.
A noticeable Reynolds number dependence has been observed
in other scalable solvers such as those discussed in~\cite{Elman2005}.
% (this behaviour for Vanka smoothers was also observed in~\cite{Vanka1986}, and
% for Braess-Sarazin in~\cite{})
Overall, the Braess-Sarazin smoother requires more iterations than the Vanka method
as in the Stokes example.
However, the gap between the two smoothers is somewhat larger than on the Stokes example,
%the s' behaviour is similar to the lid cavity
%problem example. However, the gap in the number of iterations between the two multigrid methods grows
%with increased Reynolds number
which perhaps suggests larger sensitivity of the
Braess-Sarazin smoother to higher contributions of the convective terms.
It should be noted that limitations of the Braess-Sarazin method have lead to
more modern physics-based approaches such as the pressure convection-diffusion
preconditioner and the least-squares commutator preconditioner that better account
for the convective term~\cite{Elman2005}. We hope to explore this class of physics-based
methods in the context of multigrid smoothers (as opposed to as preconditioners).

\begin{table}[t]
\caption{Number of iterations for the Navier-Stokes lid cavity.}
\label{t:ns_lidcavity}
\centering
\tabsize
\begin{tabular}{ccccccc}
  \toprule
  \multirow{2}{*}{Dofs} &
  \multicolumn{3}{c}{$\nu = 0.01$} & \multicolumn{3}{c}{$\nu = 0.002$} \\
  \cmidrule(lr){2-4} \cmidrule(lr){5-7}
  & Vanka (1,1) & BS (2,2) & ILU(1) (1,1) & Vanka (1,1) &
  BS (2,2) & ILU(1) (1,1)  \\
  \midrule
     659 & 20 & 19 &  6 & 38 & 44 &  8 \\
    2467 & 26 & 22 &  7 & 49 & 55 &  9 \\
    9539 & 22 & 26 &  9 & 53 & 60 & 13 \\
   37507 & 24 & 36 & 18 & 59 & 68 & 31 \\
  148739 & 32 & 46 & 22 & 48 & 78 & 23 \\
  \bottomrule
\end{tabular}
\end{table}

\begin{table}[t]
\caption{Number of iterations for the Navier-Stokes backward-facing step.}
\label{t:ns_backward_step}
\centering
\tabsize
\begin{tabular}{cccccccc}
  \toprule
  \multicolumn{4}{c}{$\nu = 0.02$} & \multicolumn{4}{c}{$\nu = 0.005$} \\
  \cmidrule(lr){1-4} \cmidrule(lr){5-8}
  Dofs & Vanka (1,1) & BS (2,2) & ILU(1) (1,1) & Dofs & Vanka (1,1) &
  BS (2,2) & ILU(1) (1,1)  \\
  \midrule
     479 & 20 & 21 &  6 &    889 & 51 & 56 &  5\footnotemark \\
    1747 & 22 & 26 &  8 &   3287 & 49 & 58 &  6\addtocounter{footnote}{-1}\footnotemark \\
    6659 & 25 & 24 & 11 &  12619 & 43 & 58 & 18  \\
   25987 & 26 & 32 & 17 &  49427 & 37 & 53 & 24  \\
  102659 & 26 & 47 & 25 & 195619 & 47 & 68 & 46  \\
  \bottomrule
\end{tabular}
\end{table}
\footnotetext{ILU(1) + RCM did not converge in 100 iterations; instead, ILU(2) + RCM is run}

\begin{table}[t]
\caption{Number of iterations for the Navier-Stokes obstacle.}
\label{t:ns_obstacle}
\centering
\tabsize
\begin{tabular}{ccccccc}
  \toprule
  \multirow{2}{*}{Dofs} &
  \multicolumn{3}{c}{$\nu = 0.02$} & \multicolumn{3}{c}{$\nu = 0.005$} \\
  \cmidrule(lr){2-4} \cmidrule(lr){5-7}
  & Vanka (1,1) & BS (2,2) & ILU(1) (1,1) & Vanka (1,1) &
  BS (2,2) & ILU(1) (1,1)  \\
  \midrule
     660 & 23 & 28 &  8 & 70 & 63 & 29\addtocounter{footnote}{-1}\footnotemark \\
    2488 & 23 & 33 &  8 & 41 & 62 & 17 \\
    9512 & 26 & 28 & 12 & 41 & 65 & 19 \\
   37168 & 24 & 32 & 14 & 34 & 51 & 22 \\
  146912 & 36 & 50 & 24 & 30 & 56 & 30 \\
  \bottomrule
\end{tabular}
\end{table}

\section{Conclusion}\label{s:conclusion}
A new AMG coarsening approach has been proposed for \qq{} mixed discretizations of
Stokes and Navier-Stokes equations. The advantage of the new method is in its
preservation of the spatial location relationship between pressure and velocity
unknowns throughout multigrid hierarchy, so that the qualitative structure of
the finest level is preserved on coarse levels. This is achieved by utilizing
information gathered during pressure coarsening to guide the construction
of the velocity grid transfer. The determination of grid transfer coefficients 
is then obtained by utilizing a flexible EMIN-AMG framework. A key feature of the 
proposed approach is that it coarsens fairly aggressively. In this way, the
resulting multigrid operator complexity is quite low. This implies that the 
growth in the storage and in the V cycle cost remains quite modest as
the number of multigrid levels increases.  Experiments have been 
conducted with three different smoothers to demonstrate the suitability of 
the AMG hierarchy generated by the proposed procedures. Though there is some
variation in iteration counts needed for convergence, the overall iteration/convergence rate
trends are well-behaved. In several cases, the required number of iterations does
not increase as the mesh is refined while in some other cases there is some modest
iteration growth.  Additionally, two of considered smoothers (Vanka and 
Braess-Sarazin) specifically target incompressible flow problems and rely somewhat
on sub-matrices capturing basic attributes of corresponding PDE operators, which then 
must be maintained on coarse levels of the multigrid hierarchy. 

%The proposed approach is first of its kind (to authors knowledge). 
The paper concentrated on the \qq{} approximation due to it simplicity. However, it is
hoped that %eems
%reasonable to expect 
the method can be extended to other mixed discretization methods based
on the key idea that the coarsening of one type of variable guide the coarsening
of other variables. 
%the construction of the second prolongator seems sufficiently general. The main
%effort of extending the framework would be the determination of the proper
%heuristics for the process.
One possible
%A somewhat different natural 
direction that we intend to explore are %would be a potential formulation of a
resistive magnetohydrodynamics (MHD) systems which also include electro-magnetic
effects.
%with three sets of variables:
%pressure, velocities and magnetic field. However, it should still be possible to
%use the described approach as there is no coupling between velocities and
%magnetic field variables, and thus the pressure prolongator construction can
%guide both velocity and magnetic field prolongator independently. Both
%Braess-Sarazin and Vanka smoothers are also applicable for this problem.

%\section*{Acknowledgments}
%\input{ack}
%
%\appendix
%\section{Appendix}\label{s:appendix}
%\input{appendix}

\bibliographystyle{wileyj}
\bibliography{q2q1}

\end{document}